\documentclass[12pt]{article}

\usepackage{amsmath,amsthm}
\usepackage{amssymb}
\usepackage{epsfig}

\usepackage[utf8x]{inputenc}
\usepackage{easybmat}
\usepackage{graphicx,hyperref}

\usepackage{url}

\newtheorem{theorem}{Theorem}[section]
\newtheorem{lemma}[theorem]{Lemma}
\newtheorem{prop}[theorem]{Proposition}
\newtheorem{corr}[theorem]{Corollary}

\theoremstyle{definition}
\newtheorem{define}[theorem]{Definition}
\newtheorem{remark}[theorem]{Remark}

\newcommand{\dsty}{\displaystyle}
\newcommand{\casefrac}[2]{ {\textstyle{ \frac{#1}{#2} } } }

\newcommand{\C}{ {\cal{C}} }
\newcommand{\N}{ {\cal{N}} }
\newcommand{\M}{ {\cal{M}} }
\newcommand{\Nvex}{ {\cal{N}}_{\mbox{cvex}} }
\newcommand{\Ncav}{ {\cal{N}}_{\mbox{cave}} }
\newcommand{\Cv}{ {\cal{CV}} }

\newcommand{\wl}{\hat{\lambda}}
\newcommand{\wx}{\hat{x}}
\newcommand{\wy}{\hat{y}}
\newcommand{\wa}{\hat{a}}
\newcommand{\wb}{\hat{b}}
\newcommand{\wc}{\hat{c}}
\newcommand{\wM}{\widehat{M}}

\newcommand{\pf}{\noindent {\bf Proof: }}
\newcommand{\enpf}{\hfill $\Box$ \vspace{.2in} }
\newcommand{\ts}{\vspace{0.1in}}
\newcommand{\vs}{\vspace{0.2in}}
\newcommand{\nin}{\noindent}

\begin{document}

\title{On Kite Central Configurations}

\author{Gareth E. Roberts\thanks{
Dept. of Mathematics and Computer Science,
College of the Holy Cross, groberts@holycross.edu}}

\maketitle

\begin{abstract}
We study kite central configurations in the Newtonian four-body problem.  We present a new proof that there exists a unique convex kite central configuration for a given choice of positive masses and a particular ordering of the bodies.  Our proof uses tools from differential topology (e.g., the Poincar\'{e}-Hopf Index Theorem) and computational algebraic
geometry (e.g., Gr\"{o}bner bases). We also discuss concave kite central configurations, including degenerate examples and bifurcations. 
Finally, we numerically explore the linear stability of the corresponding kite relative equilibria,
finding that the heaviest body must be at least 25 times larger than the combined masses of the other three bodies in order to be linearly stable.  
\end{abstract}



\vs

{\bf Dedication:}  In memory of Dieter S. Schmidt and his important contributions to celestial mechanics.

\section{Introduction}

Central configurations play a key role in understanding solutions of the Newtonian $n$-body problem.  From rest, a central configuration (c.c.) will collapse in on itself homothetically.
When given the correct initial velocity, a planar c.c.~rotates rigidly about its center of mass, generating a relative equilibrium periodic solution.  
Any group of bodies heading toward a collision must asymptotically approach a central configuration.  Bifurcations in the topology
of the integral manifolds in the planar problem occur at values where central configurations exist.  
For these and many other reasons, the study of central configurations is an active subfield of celestial mechanics.

In this paper we focus on four-body {\em kite} central configurations, where two of the bodies lie on an axis of symmetry and the other two
bodies (necessarily of equal mass) are positioned equidistant from this axis.  In general, a four-body planar c.c.~is called {\em convex} if 
no body lies within the convex hull of the other three.  If the configuration is not convex, it is called {\em concave}.
Kites may either be convex or concave and there are notable differences between the two cases.

There are many related results on four-body central configurations, beginning with the work of MacMillan and Bartky, who proved that for any four masses and
a particular ordering of the bodies, there exists a convex central configuration~\cite{MacB}.  Xia later gave a simpler proof of this fact~\cite{xia}.  
If a convex c.c.~has two opposite and equal masses, then it must be a kite~\cite{albouy}.  This is not true in the concave case.  
If there are two pairs of equal masses located at two adjacent vertices, then the configuration
must be an isosceles trapezoid~\cite{fern}.  In his doctoral thesis, Hampton proved that four-body concave c.c.'s exist for any choice of masses~\cite{hampton}.

One of the major open questions concerning four-body central configurations is if uniqueness holds for convex central configurations.  This is Problem \#10 on a 
well-known list of open problems~\cite{albouy-pblms}.  Using resultants and the method of rational parametrization, Leandro showed that there exists a unique 
convex kite central configuration for any choice of masses~\cite{Leandro} (assuming the masses equidistant from the axis of symmetry are equal).  
Recently, Santoprete has shown that for those masses where a co-circular central configuration
exists, such a configuration is unique~\cite{manuele1}.  He also proved this fact for trapezoidal c.c.'s~\cite{manuele2}.  
In the case where three of the masses are sufficiently small, Corbera et al.~have shown that there is a unique convex c.c.~for each ordering~\cite{CCLM}.
Sun, Xie, and You gave a numerical computer-assisted proof using interval
arithmetic and the Krawczyk operator for the uniqueness of four-body convex c.c.'s on a large subset of masses bounded away from zero~\cite{zhifu}.

Inspired by Santoprete's recent work, we give a new proof for the uniqueness of the four-body convex kite central configurations.  
We are hopeful that our proof technique will generalize to the uniqueness problem for general four-body convex c.c.'s.\footnote{Research in progress.}
While Santoprete applied Morse theory to obtain his results, we make use of the Poincar\'{e}-Hopf Index Theorem.
We follow Moeckel's approach (see~\cite{rick-book,rick-notes,rick2}) and treat kite central configurations as rest points of a certain gradient flow 
determined by the restriction of the potential function~$U$ to a special configuration space.

Here is an outline of our proof:

\begin{enumerate}
\item  Working in $\mathbb{R}^4$, we define a normalized configuration space $\N$ that is the surface of an ellipsoid.  It is important to require that the
center of mass be located at the origin.

\item  Kite central configurations are shown to be rest points of the gradient flow of $U|_{\N}$.

\item  Restricting to the set of convex configurations $\Nvex$, we show that there exists a smooth manifold $\M \subset \Nvex$ with boundary for which the gradient vector field exists and points
outward on the boundary.  This gives a simple proof for the existence of convex kite c.c.'s.

\item  We prove that any critical point in $\M$ has Poincar\'{e} index equal to one by  
showing that the product of the nontrivial eigenvalues of the modified Hessian is always positive.

\item  Finally, using the Poincar\'{e}-Hopf Index Theorem, we conclude that there can only be one critical point, thereby verifying uniqueness.
\end{enumerate}

Steps 1, 2, and 3 are accomplished in Sections 2, 3, and 4, respectively.  An important formula for the product of the non-trivial eigenvalues is derived in 
Section~\ref{subsec:Prod-evals}.  The last two steps of the proof are carried out in Section 5.  While most of our computations were done by hand, 
Gr\"{o}bner bases and symmetric variables were used for the index calculation.  These symbolic computations (no numerical estimates were necessary) 
were performed with SageMath~\cite{sage}.

In Section 6 we study the concave case and the mass map from the configuration space to the set of normalized masses.  
In contrast to the convex setting, we find a curve of degenerate central configurations and a two-to-one mass map (non-uniqueness).  These results concur with those
of Leandro~\cite{Leandro}.  The index is computed for
each type of concave configuration and the degenerate mass value for the well-studied $1+3$-gon c.c.~is recovered.

Finally, in Section 7 we carefully investigate the linear stability of the relative equilibria for both the convex and concave kites.  Although none of the
concave examples are linearly stable, there are interesting bifurcation curves that occur when the eigenvalue structure changes.  For the convex case,
we find that linearly stable examples only occur if one mass (with our assumptions it is $m_1$) is substantially larger than the others.
This agrees with Moeckel's dominant mass
conjecture (see Problem \#15 in~\cite{albouy-pblms}) as well as earlier numerical work done by Sim\'{o}~\cite{simo}.  
Specifically, we compute that the infimum of the mass ratio
$m_1/(m_2 + m_3 + m_4)$ over the set of linearly stable solutions is $(25 + 3\sqrt{69}\,)/2 \approx 24.9599$.

Figure~\ref{Fig:setup} was created with the open-source software SageMath~\cite{sage}.  All other figures were produced
using MATLAB~\cite{matlab}.

\section{Set Up}

A central configuration is a special arrangement of bodies such that the force due to gravity on each body points toward the
center of mass, with a common proportionality factor.  Let $m_0 = m_1 + \cdots + m_n$ represent the total mass.
The planar configuration $q = (q_1, q_2, \ldots, q_n)$ with $q_i \in \mathbb{R}^2$ is a {\em central configuration} if
\begin{equation}
\sum_{j \neq i}^n  \frac{m_i m_j (q_j - q_i)}{r_{ij}^3}
 \; + \; \lambda  \, m_i (q_i - c_0)  \; = \; 0
\label{cc:maineq}
\end{equation}
for some scalar $\lambda$ independent of~$i$, where $r_{ij} = || q_j - q_i ||$ is the Euclidean distance between the $i$th and $j$th bodies and 
$c_0 = \frac{1}{m_0} \sum_i m_i q_i$ is the center of mass.

One approach to the study of central configurations, which will be important for this work, is to regard c.c.'s as critical points of the Newtonian potential function
\begin{equation}
U(q)  \; = \;  \sum_{i < j}^n \; \frac{m_i m_j}{r_{ij}}
\label{eq:NewtPotential}
\end{equation}
subject to the constraint $I(q) = 1$, where 
$I$ is the moment of inertia with respect to the center of mass,  
$$
I(q) \; = \;   \sum_{i = 1}^n m_i ||q_i - c_0||^2 \, .
$$
Using these two functions, equation~(\ref{cc:maineq}) can be written more compactly as
\begin{equation}
\nabla U(q) + \casefrac{\lambda}{2} \nabla I(q) \; = \; 0,
\label{eq:cc-UI}
\end{equation}
which verifies that central configurations are critical points of $U$ restricted to the mass ellipsoid $I = 1$, with Lagrange multiplier $\lambda/2$.
Since $U$ and $I$ are homogeneous functions of degree $-1$ and $2$, respectively, equation~(\ref{eq:cc-UI}) gives a useful formula for $\lambda$, namely
\begin{equation}
\lambda  \; = \;  \frac{U(q)}{I(q)}  \, .
\label{eq:lambda-UI}
\end{equation}
Due to the symmetries of the $n$-body problem, any rotation or scaling of a central configuration generates another c.c.  Thus we think of 
central configurations as belonging to equivalence classes.  Uniqueness in the context of this article means a unique equivalence class.

We will set the problem up in $\mathbb{R}^4$, even though the defining equations can be reduced to two dimensions (see~\cite{CCR}).
After several attempts at this problem with various choices of variables, we opted to require the center of mass be at the origin and to work in $\mathbb{R}^4$.
This produced the easiest set of equations and derivatives to analyze (for example, it gives a simple formula for the moment of inertia).

Without loss of generality, assume that bodies 1 and~3 are located on the axis of symmetry of the kite, with body~1 lying 
to the right of body~3.  We will also assume that body 2 lies above the axis of symmetry.  This gives a sequential ordering of the bodies in
the convex case (counterclockwise direction).
Introduce the mass parameter $m$ with $0 < m < 1$ and let $m_2 = m_4 = m/2$.  
We normalize the masses by assuming that $m_0 = 1$ so that $m_1 + m_3 + m = 1$.

\begin{figure}[t]
\begin{center}
\includegraphics[height = 213bp]{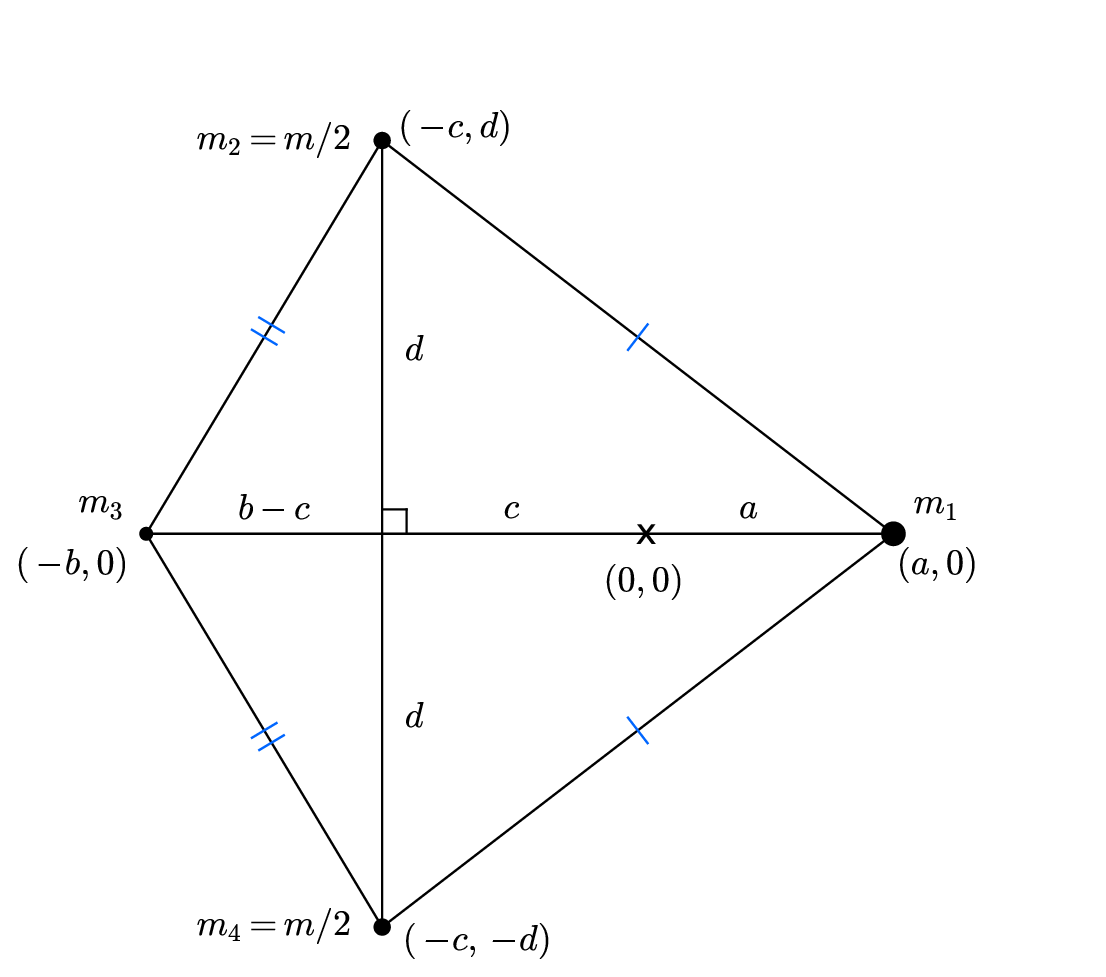}
\hspace{-0.35in}
\includegraphics[height = 210bp]{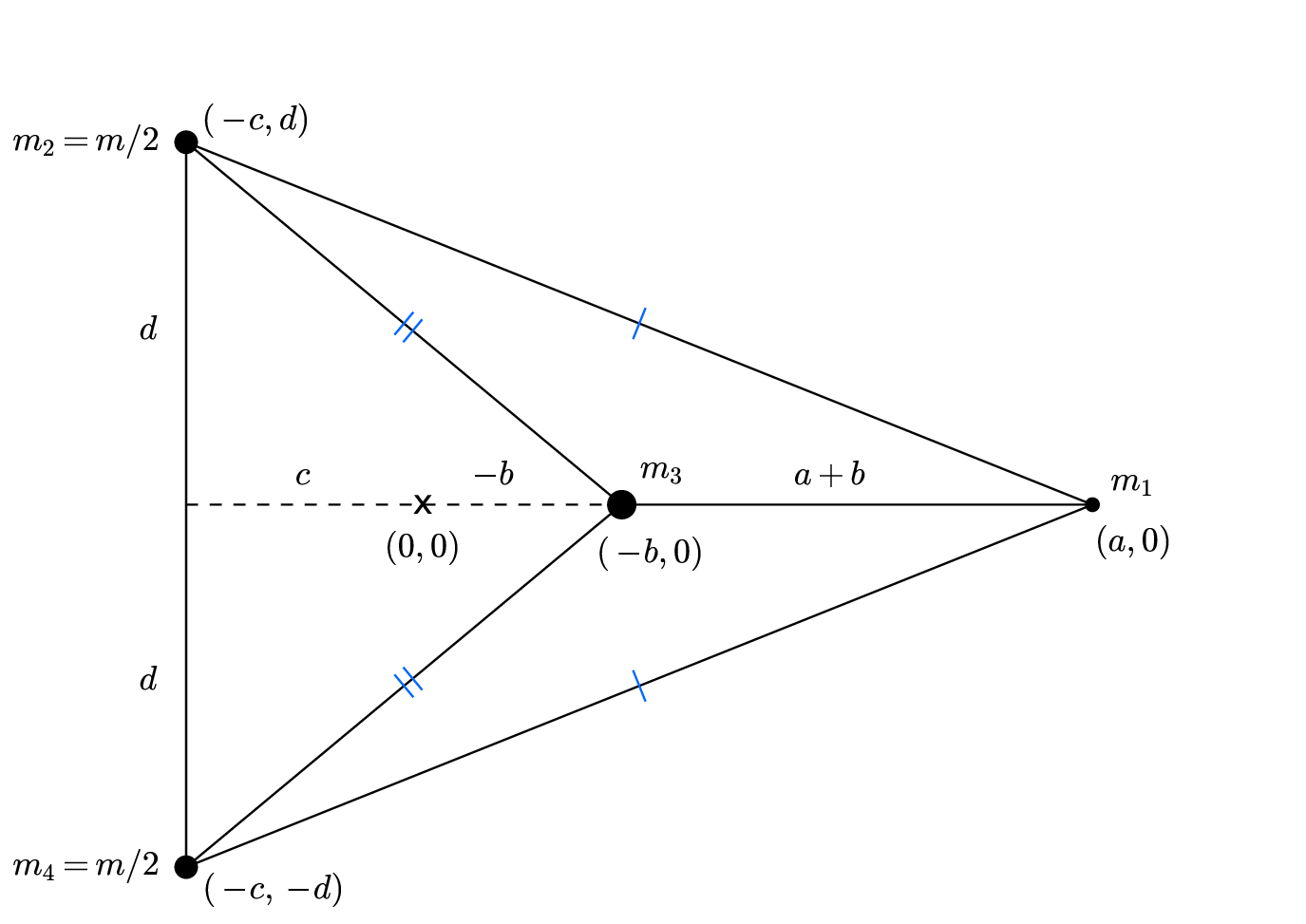}
\end{center}
\vspace{-0.2in}
\caption{Two kite configurations: the convex configuration (left) has $x = a + c > 0$ and $y = b - c > 0$, 
while the concave configuration (right) has $x = a + c > 0$ and $y = b - c < 0$. The center of mass in each figure is at the origin.}
\label{Fig:setup}
\end{figure}

Throughout the paper we let $z = (a,b,c,d) \in \mathbb{R}^4$.  The four bodies of the kite have 
positions $q_1 = (a,0), q_2 = (-c,d), q_3 = (-b,0),$ and $q_4 = (-c,-d)$ (see Figure~\ref{Fig:setup}).
Unless otherwise stated, we will assume that the quantities $a$, $a + b$, and $d$ are strictly positive.

The center of mass is at the origin provided the function
$$
C(z) \;= \;  m_1 a - m_3 b - m c \; = \;  0  \, .
$$
Instead of using this equation to eliminate one variable, we will regard $C = 0$ as a constraint on the configuration space.  Assuming that
$C = 0$, the moment of inertia reduces nicely to
$$
\widehat{I}(z) \; = \; m_1 a^2 + m_3 b^2 + mc^2 + md^2 
$$
and thus $\widehat{I} = 1$ defines an ellipsoid in $\mathbb{R}^4$.

Let $M$ be the diagonal mass matrix $M = \mbox{diag} \{m_1, m_3, m, m \}$.  Define a metric (or mass inner product) on $\widehat{I} = 1$ by
$$
<\!v, w\!> \;  =  \;  v^{\footnotesize{T}} M w, \quad \mbox{for } v,w \in \mathbb{R}^4.
$$
Denote the vector $\xi = (1, -1, -1, 0)$ and note that the center of mass condition can be written as $<z, \xi> \, = 0$.  
Then, the {\em normalized configuration space} is
$$
\N = \{ z =(a,b,c,d) \in \mathbb{R}^4: \widehat{I}(z) = \,  <\!z, z\!> \, = 1 \mbox{ and } C(z) = \, <\!z, \xi \!> \, = 0 \} .
$$
Since $C = 0$ describes a plane containing the origin, $\N$ is the surface of an ellipsoid (topologically a 2-sphere).

With respect to the mass inner product, $z$ is orthogonal to $\N$ since $\nabla \widehat{I}(z) = 2Mz$.  Likewise, $\xi$ is also orthogonal
to $\N$ since $\nabla C(z) = M\xi$. Thus, the tangent space to $\N$ at the point $z$ is simply
$$
T_z \, \N \; = \;  \{ v \in \mathbb{R}^4:  \, <\!v, z\!> \, = 0 \mbox{ and }  <\!v, \xi \!> \, = 0 \} .
$$

Kite central configurations with our prescribed ordering are critical points of $U|_{\N}$.
To make this more precise, define the linear function $\Phi: \mathbb{R}^4 \mapsto \mathbb{R}^8$ by
$$
\Phi(z) \; = \; (a,0,-c,d,-b,0,-c,-d) \, .
$$ 
Note that $\Phi$ simply maps $z$ to the positions of our kite configuration in $\mathbb{R}^8$.
Let $\widehat{U}(z) = U(\Phi(z))$ denote the potential as a function of $z \in \mathbb{R}^4$.  Assuming that $a+b > 0$ and $d > 0$, 
the mutual distances in the $z$ coordinates are given by
\begin{eqnarray}
r_{12}^2 \; = \;  (a+c)^2 + d^2   & \qquad &  r_{13} \; = \;  a + b,  \label{eq:md1} \\[0.05in]
r_{23}^2 \; = \;  (b-c)^2 + d^2   & \qquad &   r_{24} \; = \;  2d, \label{eq:md2} \\[0.05in]
r_{14} \; = \;  r_{12} & \qquad &  r_{34} \; = \;  r_{23} \,.  \nonumber
\end{eqnarray}
As is the case with general central configurations, the center of mass constraint is automatically satisfied for
critical points of $\widehat{U}$ restricted to $\widehat{I} = 1$~\cite{rick-book}.

\begin{lemma}
If $z$ is a critical point of $\widehat{U}$ restricted to $\widehat{I} = 1$, then $z \in \N$.
\label{lemma:cc-centerofmass}
\end{lemma}

\pf
This follows from the translational invariance of $U$ and $\widehat{U}$.  Note that $\Phi(z + t \xi)$ shifts each body to the right by $t$ units.  It follows that
$$
\widehat{U}(z + t \xi) \;  = \;  U(\Phi(z + t\xi)) \; = \;  U(\Phi(z)) \; = \;  \widehat{U}(z)
$$ 
for any $t \in \mathbb{R}$.  Differentiating the outer functions of this equation with respect to $t$ 
and evaluating at $t=0$ gives $\nabla \widehat{U}(z) \cdot \xi = 0$, where ``$\, \cdot \,$''
represents the standard Euclidean dot product on $\mathbb{R}^4$.

Suppose that $z$ is a critical point of $\widehat{U}$ restricted to $\widehat{I} = 1$.  Then $\nabla \widehat{U}(z) + \frac{\lambda}{2} \nabla \widehat{I}(z) = 0$
for some nonzero constant $\lambda$.  Taking the standard dot product of this equation with $\xi$ shows that $\lambda M z \cdot \xi = 0$, which
in turn implies that $C(z) = 0$.  It follows that $z \in \N$.
\enpf

Let $V = \widehat{U}|_{\N}$ be the restriction of the potential function to $\N$.
Lemma~\ref{lemma:cc-centerofmass} allows us to ignore the center of mass constraint when searching for critical points of $V$.

Define $\widehat{L}(z) = \widehat{U} + \frac{\lambda}{2} \widehat{I}$ for $z \in \mathbb{R}^4$.  If $z$ is a critical point of $\widehat{L}$, then $z$ is a 
critical point of $V$ and $z \in {\cal{N}}$ by Lemma~\ref{lemma:cc-centerofmass}.  
Computing $\nabla \widehat{L}(z) = 0$ yields the following four equations for a critical point $z$:
\begin{eqnarray}
m(a+c) r_{12}^{-3} + m_3(a + b) r_{13}^{-3} & = &  \lambda a \, ,   \label{eq:cc1} \\
m(b-c) r_{23}^{-3} + m_1(a + b) r_{13}^{-3} & = &  \lambda b   \, , \label{eq:cc2} \\
m_1(a+c) r_{12}^{-3} - m_3(b - c) r_{23}^{-3} & = &  \lambda c \, ,    \label{eq:cc3}
\end{eqnarray}
and
\begin{equation}
\lambda \; = \;  m_1 r_{12}^{-3}  + m_3 r_{23}^{-3} + m r_{24}^{-3}  \, .  \label{eq:lambda}
\end{equation}
Here, the mutual distances are given by formulas (\ref{eq:md1}) and~(\ref{eq:md2}).  
Equation~(\ref{eq:cc3}) is superfluous, since it follows directly from equations (\ref{eq:cc1}), (\ref{eq:cc2}), and $C(z) = 0$, but it is useful for
proving the following lemma.

\pagebreak

\begin{lemma} 
$z$ is a critical point of $V$ if and only if $\Phi(z)$ is a kite central configuration.
\label{lemma:cc}
\end{lemma}

\pf
Suppose that $\Phi(z)$ is a central configuration.  By equations (\ref{eq:cc-UI}) and~(\ref{eq:lambda-UI}), $\Phi(z)$ is a critical point of 
$L(q) = U + \frac{\lambda}{2} I, q \in \mathbb{R}^8$, where $\lambda = U(\Phi(z))/I(\Phi(z))$.
By Lemma~\ref{lemma:cc-centerofmass}, critical points of $\widehat{L} = \widehat{U} + \frac{\lambda}{2} \widehat{I}$ are automatically critical points of $V$.

We have $\widehat{L}(z) = L(\Phi(z))$ by definition. By the chain rule, this implies that
$
\nabla \widehat{L}(z) \; = \;  (D\Phi)^{\footnotesize{T}} \nabla L( \Phi(z) )  ,
$
where
$$
 (D\Phi)^{\footnotesize{T}} \; = \;  
\begin{bmatrix}
1 & 0 & 0 & 0     & 0 & 0 & 0 & 0 \\
0 & 0 & 0 & 0     & -1 & 0 & 0 & 0 \\
0 & 0 & -1 & 0     & 0 & 0 & -1 & 0 \\
0 & 0 & 0 & 1     & 0 & 0 & 0 & -1 
\end{bmatrix}
.
$$
Since $\nabla L( \Phi(z) ) = 0$, it follows immediately that $z$ is a critical point of $\widehat{L}$, as desired.

Conversely, suppose that $z$ is a critical point of $V$.  Then $z$ is a critical point of $\widehat{L} = \widehat{U} + \frac{\lambda}{2} \widehat{I}$, where
$\lambda = \widehat{U}(z)/\widehat{I}(z) = \widehat{U}(z)$, $z \in {\cal{N}}$, and equations (\ref{eq:cc1}) through~(\ref{eq:lambda}) are satisfied. 
Defining $L = U + \frac{\lambda}{2}I$ as above, we see that $\nabla L( \Phi(z) )$ is in the
kernel of $(D\Phi)^{\footnotesize{T}}$.  This means that $\nabla L( \Phi(z) )$ has the form $(0, l_2, l_3, l_4, 0, l_6, -l_3, l_4)$. We want to show that 
the $l_i$'s are all zero.

Write $q = (q_1, q_2, q_3, q_4) \in \mathbb{R}^8$.  By symmetry, it is clear that the second components of $\partial L/\partial q_1(\Phi(z))$ and 
$\partial L/\partial q_3(\Phi(z))$ vanish.  This shows that $l_2 = l_6 = 0$.  The vector $\partial L/\partial q_2(\Phi(z))$ also vanishes by  
equations (\ref{eq:cc3}) and~(\ref{eq:lambda}). This shows that $l_3 = l_4 = 0$ and thus $\nabla L( \Phi(z) ) = 0$.  Consequently,
$\Phi(z)$ is a central configuration.  The fact that it is a kite follows since $r_{14} = r_{12}$ and $r_{23} = r_{34}$ hold for any configuration $\Phi(z)$.
\enpf

\section{Gradient Flow and the Modified Hessian}

Based on Lemma~\ref{lemma:cc}, we will refer to critical points of $V$ as kite central configurations with the understanding that $z \in \mathbb{R}^4$ is being
identified with the full configuration $\Phi(z) \in \mathbb{R}^8$.  We regard critical points of $V$ as rest points of the gradient flow on the normalized
configuration space ${\cal{N}}$.  This formulation will be particularly useful in the convex setting.

For ease of notation we drop the $\, \widehat{ \, } \, $ from the reduced potential and inertia functions.
As a reminder, for $z = (a,b,c,d) \in \mathbb{R}^4$, we have
$$
U(z) \; = \;  \frac{m m_1}{r_{12}} +  \frac{m m_3}{r_{23}} +  \frac{m_1 m_3}{r_{13}} + \frac{m^2}{4 r_{24}}  \quad \mbox{ and }  \quad
I(z) \; = \;  m_1 a^2 + m_3 b^2 + mc^2 + md^2  \, ,
$$
where the mutual distances are given by equations (\ref{eq:md1}) and~(\ref{eq:md2}).

\begin{lemma}
Let $X(z)$ denote the vector field $M^{-1} \nabla U(z) + U(z) z$ for $z \in \mathbb{R}^4$.  Then
$X|_{\N}$ is the gradient vector field of $V = U|_{\N}$ with respect to the metric $<\!v, w\!> \;  =  \;  v^{\footnotesize{T}} M w$
and rest points (zeros) of $X$ are kite central configurations.
\label{lemma:gradfield}
\end{lemma}

\pf
First we check that $X$ is tangent to ${\N}$.  Recall that 
$$
T_z \, {\cal{N}} \; = \;  \{ v \in \mathbb{R}^4:  \, <\!v, z\!> \, = 0 \mbox{ and }  <\!v, \xi \!> \, = 0 \} \, .
$$
Using the homogeneity of $U$, for any $z \in {\N}$, we have
\begin{eqnarray*}
<\!X, z\!> & = &  <\!M^{-1} \nabla U(z), z \!> + U(z) <\!z, z\!> \\
& = &  (\nabla U(z))^{\footnotesize{T}} M^{-1} M z + U(z) \cdot 1 \\
& = &  -U(z) + U(z)  \; = \; 0
\end{eqnarray*}
and
\begin{eqnarray*}
<\!X, \xi\!> & = &  <\!M^{-1} \nabla U(z), \xi \!> + U(z) <\!z, \xi\!> \\
& = &  (\nabla U(z))^{\footnotesize{T}} M^{-1} M \xi + 0 \\
& = &  (\nabla U(z))^{\footnotesize{T}} \xi  \; = \; 0 \, ,
\end{eqnarray*}
where the last step follows from the translational invariance of $U$ (see the proof of Lemma~\ref{lemma:cc-centerofmass}).
This shows that $X(z) \in T_z \, \N$.

Next we show that for $v \in T_z \, \N$, the derivative of $U$ in the direction of $v$, denoted $DU(z) v$, is equivalent to $<\!X(z), v\!>$.
Suppose that $v \in T_z \, {\cal{N}}$, which implies that $<\!v, z\!> \, = 0$.  Then
\begin{eqnarray*}
<\!X, v \!> & = &  <\!M^{-1} \nabla U(z), v \!> + U(z) <\!z, v \!> \\
& = &  (\nabla U(z))^{\footnotesize{T}} M^{-1} M v + 0 \\
& = &  DU(z) v  \, .
\end{eqnarray*}

Finally, if $z \in \N$ is a rest point of $X$, then $\nabla U(z) + U(z) M z = 0$.  This is equivalent to $\nabla U(z) + \frac{\lambda}{2} \nabla I(z)$, since
$\lambda = U(z)/I(z)$ reduces to $\lambda = U(z)$ on $\N$.  Therefore, rest points of $X$ are equivalent to critical points of $V$, as desired.
\enpf

In order to determine the type of critical point, we need a matrix representation of the Hessian of~$V$ (see pp.~48--50 of~\cite{rick-notes}).
Consider the function $g(z) = \sqrt{I(z)\,} \, U(z)$ on $\mathbb{R}^4$.  
Note that $g|_{\N} = U|_{\N} = V$.  If $z \in \N$ (which means $I(z) = 1$) and $z$ is a central configuration, then we compute
$$
Dg(z) v \; = \;  (\nabla U(z))^{\footnotesize{T}} v + U(z) (M z)^{\footnotesize{T}} v \; = \;  <\!X(z), v \!>  \; = \; 0 \, .
$$
Thus critical points of $g$ are c.c.'s.  In addition, if $z \in \N$ and $z$ is a central configuration, then for any vector $v \in T_z \, \N$, we find
\begin{eqnarray*}
D^2g(z) v & = & \left[  (\nabla U(z) - U(z)Mz) (Mz)^{\footnotesize{T}} +  U(z) M  + (Mz) (\nabla U(z))^{\footnotesize{T}} + D^2U(z)  \right] v \\[0.07in]
& = &  \left[  (\nabla U(z) - U(z)Mz) z^{\footnotesize{T}}M + U(z) M +  Mz (-U(z) Mz)^{\footnotesize{T}} + D^2U(z)    \right] v \\[0.07in]
& = &   (\nabla U(z) - U(z)Mz)  <\!z, v\!>  + U(z) Mv - U(z)Mz <\!z, v\!> + D^2U(z)v  \\[0.07in]
& = &  (D^2U(z) + U(z) M)v  \, .
\end{eqnarray*}
This justifies the following definition.

\begin{define}
The {\em Hessian} of $V = U|_{\cal{N}}$ at a critical point $z$ is the restriction of the $4 \times 4$ matrix 
$$
H(z) \; = \;  D^2 U(z) + \lambda M, \quad \mbox{where } \lambda = U(z),
$$
to $T_z \, \N.$  The {\em modified Hessian} matrix is 
$$
M^{-1} H(z) \; = \;  M^{-1} D^2U(z) + \lambda I_4 \, ,
$$ 
where $I_4$ is the $4 \times 4$ identity matrix.
\end{define}

\ts

The matrix $M^{-1}H(z)$ is symmetric with respect to the mass inner product and thus has all real eigenvalues.  
Using Sylvester's Law of Inertia~\cite{meyer-book}, the number of eigenvalues with a particular sign ($+, -,$ or $0$) for the Hessian and the modified Hessian are identical.
A proof of this fact is given in Lemma 2.3 in~\cite{gr-morsethy}.  In particular, $z$ is a local minimum of $V$ 
if and only if the modified Hessian is positive definite (all eigenvalues are positive).

\subsection{Trivial eigenvalues and the derivative of $X(z)$}

There are two advantages of working with the modified Hessian.  First, it always contains two ``trivial'' positive eigenvalues, $\lambda$ and $3 \lambda$.
Second, when restricted to the tangent space $T_z \, \N$, it is equivalent to the derivative of the vector field $X$ evaluated at a rest point $z$.
We verify both of these facts below.

Due to the translational invariance of $U$, we have $D^2U(z) \xi = 0$ and consequently $M^{-1} H(z) \xi = \lambda \xi$.  This eigenvector results from the
conservation of the center of mass.  
Next, due to the homogeneity of $U$ and the fact that $z$ is a central configuration, we have that $z$ is also
an eigenvector of the modified Hessian, with eigenvalue $3\lambda$.
To see this, differentiate the identity $\nabla U(z) \cdot z = -U(z)$ with respect to~$z$ (here we treat $z$ as a variable).  We compute
\begin{equation}
D^2U(z) z + \nabla U(z) \; = \;  - \nabla U(z)  \quad \Longrightarrow \quad  D^2U(z) z = -2 \nabla U(z) \, .
\label{eq:D2Uidentity}
\end{equation}
If $z$ is a critical point of $V$, then $\nabla U(z) + \lambda M z = 0$.  Substituting this into identity~(\ref{eq:D2Uidentity}) yields $M^{-1} H(z) z = 3 \lambda z$.
The fact that $z$ is an eigenvector reflects the scaling symmetry of equations (\ref{eq:cc1})--(\ref{eq:lambda}).  We note the importance of
requiring the center of mass to be located at the origin; otherwise $z$ would not be an eigenvector of the modified Hessian.

\begin{define}
For any critical point $z \in \N$, the modified Hessian $M^{-1}H(z)$ always has two trivial eigenvalues $\lambda$ and $3 \lambda$, where $\lambda = U(z)$.
If the remaining two eigenvalues are nonzero, then $z$ is a {\em nondegenerate} critical point of $V = U|_{\N}$.
\label{def:CC}
\end{define}

Next we consider $dX(z)$, the derivative of the gradient vector field $X$ evaluated at a rest point~$z$.  Since $z$ is a central configuration, we know 
$\nabla U(z) = -U(z) M z$.  We compute
\begin{eqnarray*}
dX(z) & = &  M^{-1} D^2U(z) + z  (\nabla U(z))^{\footnotesize{T}} + U(z) I_4  \\[0.07in]
& = &   M^{-1} D^2U(z) + z (-U(z) M z)^{\footnotesize{T}} + U(z) I_4 \\[0.07in]
& = &   M^{-1} D^2U(z) - \lambda z z^{\footnotesize{T}} M   +  \lambda I_4 \\[0.07in]
& = &  M^{-1} H(z) - \lambda z z^{\footnotesize{T}} M \, .
\end{eqnarray*}
Recall that for $z \in \N$, we have $<\!z, \xi\!> \, = 0$ and $<\!z, z\!> \, = 1$.
It follows that $dX(z)$ also has $\xi$ and $z$ as trivial eigenvectors, with eigenvalues $\lambda$ and $2\lambda$, respectively.

Given a critical point $z$, define the two-dimensional linear subspace $B = \mbox{span}\{z, \xi \} \subset \mathbb{R}^4$ 
and its orthogonal complement with respect to $<\! \ast, \ast, \!>$,
$$
B^\perp  \; = \;  \{w \in \mathbb{R}^4: \, <\!w, v\!> \, = w^{\footnotesize{T}} M v = 0 \; \; \forall v \in B\} .
$$ 
The nontrivial eigenvectors for $M^{-1}H(z)$ and $dX(z)$ lie in $B^\perp$.

\begin{lemma}
The vector space $B^\perp$ has the following important properties:
\begin{itemize}
\item $B^\perp$ is equivalent to $T_z \, \N$.

\item $B^\perp$ is invariant under both the modified Hessian and $dX(z)$.

\item The matrices $M^{-1} H(z)$ and $dX(z)$ are identical when restricted to $B^\perp$.
\end{itemize}
\label{lemma:invBperp}
\end{lemma}

\pf
Recall that $T_z \, \N \; = \;  \{ w \in \mathbb{R}^4:  \, <\!w, z\!> \, = 0 \mbox{ and}  <\!w, \xi \!> \, = 0 \} .$  Since $z$ and $\xi$ are a basis for $B$, it follows
that $w \in B^\perp$ if and only if $w \in T_z \, \N$.

Next suppose that $w \in B^\perp$ and consider the vector $M^{-1}H(z) w$.  
Let $v$ be an arbitrary vector in $B$, so that $v = c_1 z + c_2 \xi$ for some $c_i \in \mathbb{R}$.  Then we have
\begin{eqnarray*}
(M^{-1}H(z) w)^{\footnotesize{T}} M v & = &  w^{\footnotesize{T}} H(z) M^{-1} M v \\[0.07in]
& = &  c_1 w^{\footnotesize{T}} H(z) z +  c_2 w^{\footnotesize{T}} H(z) \xi  \\[0.07in]
& = &  c_1 w^{\footnotesize{T}} (3 \lambda M z) + c_2 w^{\footnotesize{T}} (\lambda M \xi) \\[0.07in]
& = &  3 \lambda c_1 <\!w,z\!> + \, \lambda c_2 <\!w,\xi \!>  \\
& = & 0 \, ,
\end{eqnarray*}
since $w \in T_z \, \N$.  Hence, $M^{-1}H(z) w \in B^\perp$, which shows that $B^\perp$ is invariant under $M^{-1}H(z)$.

Finally, for any vector $w \in B^\perp$, we have
$$
-\lambda z z^{\footnotesize{T}} M w \; = \;  -\lambda z \left(z^{\footnotesize{T}} M w \right) \; = \;  -\lambda <\!w,z\!> z \; = \; 0 \, ,
$$
which implies that $dX(z) w = M^{-1}H(z) w$.  This shows the last item in the lemma and also verifies the invariance of $B^\perp$ under $dX(z)$.
\enpf

\begin{remark}
Since $X(z)$ is a vector field on the manifold $\N$, its derivative $dX$ can be viewed as a linear mapping from $T_z \, \N$ to itself.  Consequently, the invariance
of $B^\perp = T_z \, \N$ is to be expected.
\end{remark}

Lemma~\ref{lemma:invBperp} shows that $\mathbb{R}^4 = B \oplus T_z \, \N$ gives a decomposition of $\mathbb{R}^4$ into
invariant subspaces, where the trivial eigenvectors belong to $B$ and the nontrivial eigenvectors are contained in the tangent space.
Moreover, the two nontrivial eigenvalues of the modified Hessian are equivalent to those of $dX(z)$.  We will show that for the convex case, they are always positive.

\subsection{A formula for the product of the nontrivial eigenvalues}
\label{subsec:Prod-evals}

Suppose that $z$ is a critical point of $V$ and let $A = M^{-1}D^2U(z)$.
After some simplification, using algebra tricks such as $r_{ij}^{-3} = r_{ij}^{-5} \cdot r_{ij}^2$ and formulas (\ref{eq:md1}) and~(\ref{eq:md2}), we compute that $A$ is
\begingroup\makeatletter\def\f@size{10}\check@mathfonts
\def\maketag@@@#1{\hbox{\m@th\large\normalfont#1}}%
$$
\hspace{-0.2in}
\begin{bmatrix}
\zeta_1  & 2m_3 r_{13}^{-3} &  m r_{12}^{-5} \left( 2(a+c)^2 - d^2 \right) &  3md(a+c)r_{12}^{-5}    \\[0.07in]
2m_1 r_{13}^{-3} &  \zeta_2 & -m r_{23}^{-5} \left( 2(b-c)^2 - d^2 \right)  &  3md(b-c)r_{23}^{-5}    \\[0.07in]
m_1 r_{12}^{-5} \left( 2(a+c)^2 - d^2 \right)   & -m_3 r_{23}^{-5} \left( 2(b-c)^2 - d^2 \right)  & \zeta_3  &  3m_1 d(a+c) r_{12}^{-5} - 3m_3 d (b-c) r_{23}^{-5}   \\[0.07in]
3m_1 d(a+c) r_{12}^{-5} &  3m_3 d(b-c)r_{23}^{-5} &  3m_1 d(a+c) r_{12}^{-5} - 3m_3 d (b-c) r_{23}^{-5} & \zeta_4   
\end{bmatrix}
$$
\endgroup
where
\begin{eqnarray*}
\zeta_1 & = &  mr_{12}^{-5} ( 2(a+c)^2 - d^2) + 2m_3 r_{13}^{-3}  \\[0.05in]
\zeta_2 & = &  mr_{23}^{-5} ( 2(b-c)^2 - d^2) + 2m_1 r_{13}^{-3} \\[0.05in]
\zeta_3 & = &  m_1 r_{12}^{-5} ( 2(a+c)^2 - d^2) + m_3 r_{23}^{-5} ( 2(b-c)^2 - d^2) \\[0.05in]
\zeta_4 & = &  3m_1 d^2 r_{12}^{-5} + 3m_3 d^2 r_{23}^{-5}  + 3m r_{24}^{-3} - \lambda \, .
\end{eqnarray*}
It is straight-forward to check that $\xi$ is in the kernel of $A$ and that $Az = 2\lambda z$ when $z$ is a solution to equations (\ref{eq:cc1})--(\ref{eq:lambda}).
As expected, this confirms that 0 and $2 \lambda$ are trivial eigenvalues of $A$.

Below we give a formula for the product of the two nontrivial eigenvalues of the modified Hessian in terms of the traces of $A$ and $A^2$.
This is easier than working with a basis for the tangent space or computing the determinant of~$A$ directly.

Let $\mu_1$ and $\mu_2$ be the nontrivial eigenvalues of $A$.  Then $\mu_1 + \lambda$ and $\mu_2 + \lambda$ are the nontrivial eigenvalues of the
modified Hessian or equivalently, of $dX(z)$.
Denote tr$(\ast)$ as the trace of a matrix.  First, note that the sum of the eigenvalues of $A$ is $0 + 2\lambda + \mu_1 + \mu_2$.  This is equivalent to the trace of $A$.  
Next, let $p(\mu)$ be the characteristic polynomial of~$A$.  We have
\begin{equation}
p(\mu) \; = \;  \mu (\mu - 2\lambda)(\mu - \mu_1)(\mu - \mu_2) \, .
\label{eq:charpolyA}
\end{equation}
Using the Leverrier-Souriau-Frame algorithm (see p.~504 in~\cite{meyer-book}), the coefficient of the quadratic term of $p(\mu)$ is equal to
$\frac{1}{2} \left[ ( \, \mbox{tr}(A) \, )^2 - \mbox{tr}(A^2) \right]$.  On the other hand, we can obtain this coefficient by expanding~(\ref{eq:charpolyA}).  This implies that
$$
\mu_1 \mu_2 + 2\lambda(\mu_1 + \mu_2) \; = \;  \frac{1}{2} \left[  ( \, \mbox{tr} (A) \, )^2 - \mbox{tr}( A^2 ) \right]  \, .
$$
Finally, we compute
\begin{eqnarray}
(\mu_1 + \lambda)(\mu_2 + \lambda) & = &  \mu_1 \mu_2 + \lambda(\mu_1 + \mu_2) + \lambda^2  \nonumber  \\[0.07in]
& = &  \frac{1}{2} \left[  ( \, \mbox{tr} (A) \, )^2 - \mbox{tr}( A^2 ) \right]  - \lambda(\mu_1 + \mu_2) + \lambda^2   \nonumber \\[0.07in]
& = &  \frac{1}{2} \left[  ( \, \mbox{tr} (A) \, )^2 - \mbox{tr}( A^2 ) \right]  - \lambda( \mbox{tr}(A) - 2\lambda) + \lambda^2  \nonumber \\[0.07in]
& = &  \frac{1}{2} \left[  ( \, \mbox{tr} (A) \, )^2 - \mbox{tr}( A^2 ) \right]  - \lambda \mbox{tr}(A) + 3\lambda^2 \, .  \label{eq:prodEvals}
\end{eqnarray}

\subsection{The Dziobek equations}

Before proceeding to an analysis of the convex kite central configurations, we show how to derive the well-known Dziobek equations~\cite{dzio}
from equations (\ref{eq:cc1}), (\ref{eq:cc2}), and~(\ref{eq:lambda}).  For a four-body kite central configuration with
$r_{12} = r_{14}$ and $r_{23} = r_{34}$, there is only one Dziobek equation, given by
\begin{equation}
(\lambda - r_{12}^{-3})(\lambda - r_{23}^{-3})  \; = \; (\lambda - r_{13}^{-3})(\lambda - r_{24}^{-3}) \, .  \label{eq:Dziobek}
\end{equation}
This is typically derived by viewing the mutual distances $r_{ij}$ as independent variables and using
the Cayley-Menger determinant as a constraint (see for example \cite{HRS, schmidt}).

First, substitute the formula for $\lambda$ from equation~(\ref{eq:lambda}) into equations (\ref{eq:cc1}) and~(\ref{eq:cc2}), and eliminate $mc$
using the equation $C = 0$.  This gives the mass ratios
\begin{equation}
\frac{m_1}{m} \; = \;  \frac{b(r_{23}^{-3} - r_{24}^{-3})}{\sigma} \quad \mbox{ and} \quad  \frac{m_3}{m} \; = \;  \frac{a(r_{12}^{-3} - r_{24}^{-3})}{\sigma}
\label{eq:massratios1}
\end{equation}
where $\sigma = a(r_{23}^{-3} - r_{13}^{-3}) + b(r_{12}^{-3} - r_{13}^{-3})$.  Next, these two formulas imply that
\begin{equation}
\frac{1}{m} \; = \;  \frac{m_1}{m} + \frac{m_3}{m} + 1 \; = \;  \frac{ b(r_{23}^{-3} - r_{24}^{-3}) + a(r_{12}^{-3} - r_{24}^{-3}) + \sigma}{\sigma}  \, .
\label{eq:massratios2}
\end{equation}
Substituting formulas (\ref{eq:massratios1}) and (\ref{eq:massratios2}) into
$$
\frac{\lambda}{m} \; = \;  \frac{m_1}{m} r_{12}^{-3} + \frac{m_3}{m} r_{23}^{-3} + r_{24}^{-3} 
$$
and solving for $\lambda$ gives
\begin{equation}
\lambda \; = \;  \frac{ r_{12}^{-3} r_{23}^{-3}  - r_{13}^{-3} r_{24}^{-3} }{ r_{12}^{-3} + r_{23}^{-3} - r_{13}^{-3} - r_{24}^{-3} } \, .
\label{eq:lambda-dist}
\end{equation}
Finally, if we cross multiply, arrange the terms with $r_{12}^{-3}$ and $r_{23}^{-3}$ on one side and $r_{13}^{-3}$ and $r_{24}^{-3}$ on the other,
and add $\lambda^2$ to both sides, we obtain Dziobek's equation~(\ref{eq:Dziobek}).  This equation will be useful for calculating the mass
ratios in Section~\ref{subsec:reducedC}.

\section{Convex Kites: A Smooth Manifold Without Collisions}

In order to apply the Poincar\'{e}-Hopf Index Theorem, we need to work on a smooth manifold.
Here we show how to shrink the configuration space slightly in order to avoid collisions and smooth the boundary of the manifold of convex configurations.  In the process, we obtain
a simple proof for the existence of convex kites reminiscent of an argument due to Xia~\cite{xia}.

Recall that the normalized configuration space $\N \subset \mathbb{R}^4$ is the ellipsoid $I = 1$ lying in the hyperplane $C=0$.  
It is helpful to study $\N$ in a three-dimensional space.  
Toward that end, introduce the variables $x = a + c$ and $y = b - c$.  In these variables, the center of mass condition $C = 0$ becomes 
\begin{equation}
c \;  =  \; m_1 x - m_3 y \, 
\label{eq-c}
\end{equation}
and the moment of inertia is
\begin{equation}
I \; = \;  mm_1 x^2 + mm_3 y^2 + m_1 m_3 (x+y)^2 + m d^2 \; = \;  1 \, .
\label{eq:newI}
\end{equation}
Points on $\N$ can be described locally by the coordinates $(x, y, d)$. The corresponding point $z = (a,b,c,d) \in \mathbb{R}^4$ is given by the image of the linear map
\begin{eqnarray}
a & = &  (1-m_1)x + m_3 y  \label{eq:a} \\
b & = &  m_1 x + (1 - m_3) y \\
c & = &  m_1 x - m_3 y \\
d & = & d \, .  \label{eq:d}
\end{eqnarray}

The signs of $x$, $y,$ $x+y$, and $d$ determine the type and ordering of the corresponding kite configuration.
Convex configurations must satisfy $xy > 0$, while concave configurations satisfy $x y < 0$.  For example, if $x > 0$, $y < 0$, $x+y > 0$ and $d > 0$, then the configuration
is concave with body~3 lying inside the triangle determined by bodies 1, 2, and~4 (see Figure~\ref{Fig:setup}).  Flipping the sign of $d$ interchanges bodies
2 and~4, which alters the ordering of the bodies.

We now focus on the convex case.  Without loss of generality, we assume that the bodies are ordered consecutively $(1 \, 2 \, 3 \, 4)$ in the counterclockwise direction,
with $a > 0$, as in the left image in Figure~\ref{Fig:setup}.
Such configurations lie in the interior of 
$$
\Nvex \; = \;  \{ z \in \N: x \geq 0, y \geq 0, \mbox{and } d \geq 0 \} .
$$
$\Nvex$ is a triangular region on the surface of the ellipsoid $I = 1$. 
Its boundary consists of three curves, two of which contain configurations that lie on the border between convex and concave configurations:
\begin{itemize}
\item  $x = 0$, where bodies 1, 2, and 4 are collinear with body 1 in the center,

\item  $y = 0$, where bodies 2, 3, and 4 are collinear with body 3 in the center,

\item  $d = 0$, where bodies 2 and 4 collide ($r_{24} = 0$). 
\end{itemize}
A collision also occurs on the boundary if $x = y = 0$ ($r_{13} = 0$, so bodies 1 and~3 collide).  
$\Nvex$ has three corners at the points 
$$
P_1 = (x_m, 0 , 0), \; P_2 = (0, y_m, 0), \; \mbox{and} \; P_3 = (0, 0, d_m) \, ,
$$
where $x_m = 1/\sqrt{m_1(1-m_1)}$, 
$y_m = 1/\sqrt{m_3(1-m_3)}$, and $d_m = 1/\sqrt{m\,}$ are the maximum possible values of $x,y,$ and $d$, respectively, on $\Nvex$.
$P_3$ corresponds to a binary collision between bodies 1 and~3, while 
$P_1$ and $P_2$ each represent a triple collision.  
Since $V$ blows up at any collision, we construct a slightly smaller
configuration space for $\Nvex$ that has no corners and avoids collisions.

\begin{figure}[h]
\begin{center}
\includegraphics[height = 240bp]{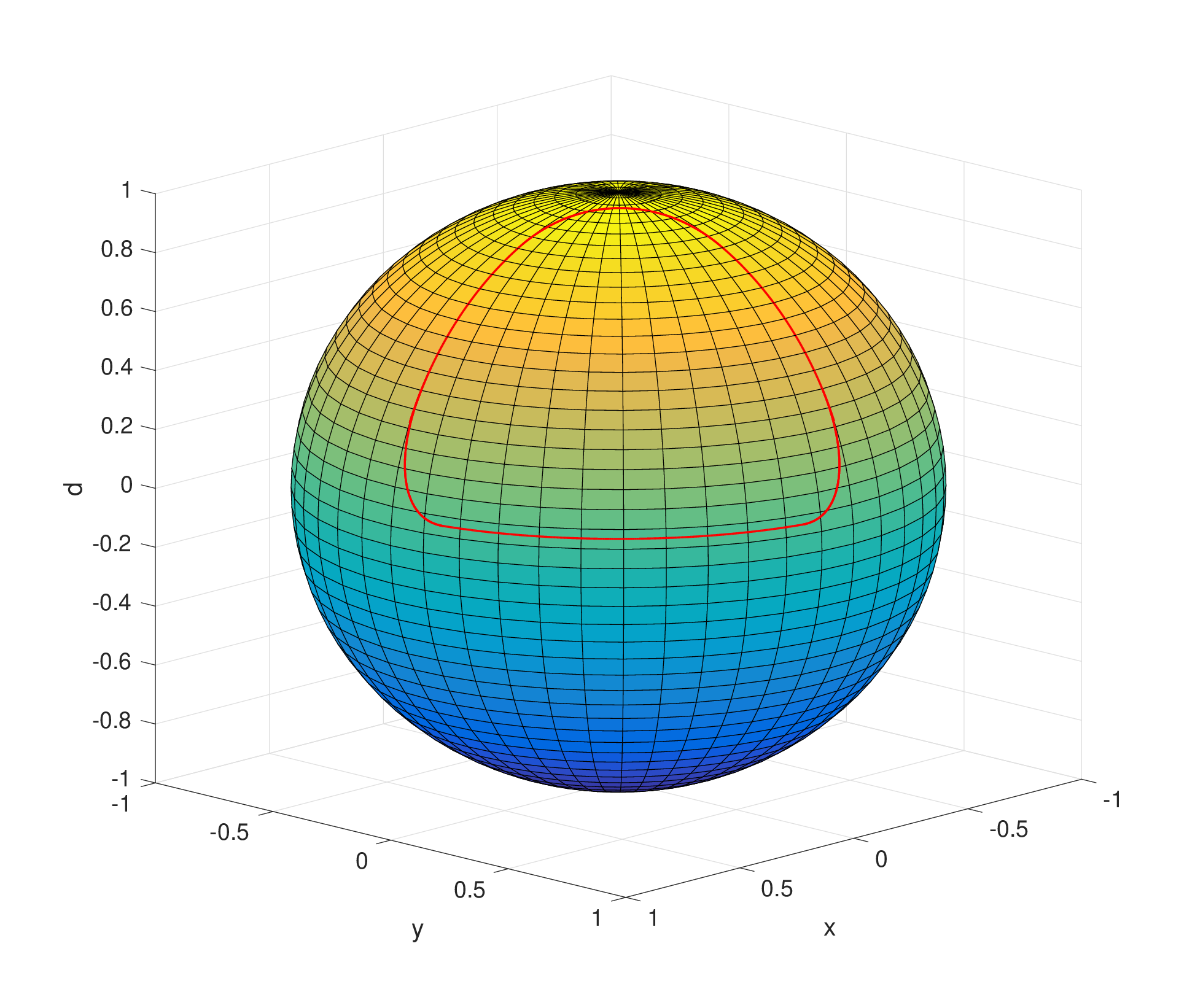}
\end{center}
\vspace{-0.3in}
\caption{An example of a smooth deformation $\gamma$ (red curve) of the boundary $\Nvex$ that avoids collisions.  The figure is shown in local coordinates $(x,y,d)$.} 
\label{Fig:SampleM}
\end{figure}

Let $\delta > 0$ be a small number.  Define $\gamma$ to be a smooth, closed curve on $\Nvex$ that contains $x=0$, $y = 0$, and $d = \delta$, and that smoothly
connects these arcs with small curves that avoid the three collisions points $P_i$.  We may take $\gamma$ as close to the boundary of $\Nvex$ as we like.
This will be the boundary of our smooth manifold. See Figure~\ref{Fig:SampleM} for an example of such a curve in local coordinates $(x,y,d)$.

In order to prove that the gradient flow points outward along $\gamma$, we first derive a general formula for the outward pointing normal in $T_z \, \N$.
Let $z = (a,b,c,d)$ be a point in $\N$.  Recall that $x = a+c$ and $y = b-c$.  The vectors $w_1$ and $w_2$ form a basis for $T_z \, \N$, where
$$
w_1 \; = \;  
\begin{bmatrix}
m m_3 y  \\ 
-m m_1 x  \\  
m_1m_3 (x+y) \\
0 
\end{bmatrix}
\quad \mbox{and} \quad
w_2 \; = \;  
\begin{bmatrix}
-m m_3 d  \\ 
-m m_1 d  \\  
0 \\
m_1 m_3 (x+y) 
\end{bmatrix} .
$$
Suppose that $\gamma$ is parametrized by $\gamma(t) = (a(t), b(t), c(t), d(t))$. Then $\gamma'(t) = \frac{d}{dt} (\gamma(t))$ lies in the tangent space and is tangent to $\gamma$.  The vector
$$
\widehat{w}(t) \; = \;  <\!\gamma'(t), w_2\!> w_1 \: - \,  <\!\gamma'(t), w_1\!> w_2
$$
is orthogonal to $\gamma'(t)$ and also lies in $T_z \, \N$.  Thus, for a fixed parameter value $t = t_0$, either $\widehat{w}(t_0)$ or $-\widehat{w}(t_0)$ is the outward pointing normal 
for the region interior to $\gamma$ at the point $z_0 = \gamma(t_0)$.
We compute that
\begin{equation}
\widehat{w}(t_0)  \; = \;  
\begin{bmatrix}
\widehat{w}_1  \\[0.07in] 
\widehat{w}_2   \\[0.07in]   
\widehat{w}_3  \\[0.07in] 
\widehat{w}_4  
\end{bmatrix} 
\; = \;  -m m_1 m_3 (x + y)  
\begin{bmatrix}
m m_3 (d \, y' - d' y)  \\[0.07in]   
- m m_1 (d \, x' - d' x)  \\[0.07in] 
m_1 m_3 \left( d (x' + y') - d' (x + y) \right) \\[0.07in] 
m_1 m_3 (y \, x' - y' x) 
\end{bmatrix} ,
\label{eq:normalVec}
\end{equation}
where all of the variables and their derivatives are evaluated at $t_0$.

Suppose that $z_0 \neq P_i$ lies on the boundary $x = 0$ of $\Nvex$.  Then we have $a(t_0) + c(t_0) = x(t_0) = 0$
and $a'(t_0) + c'(t_0) = x'(t_0) = 0$.  Using equation~(\ref{eq:normalVec}) we see that
$$
\widehat{w}_x \; = \;  -
\begin{bmatrix}
m  \\
0 \\
m_1 \\
0
\end{bmatrix} 
$$
is an outward pointing normal vector. The negative sign indicates that $x = a + c$ will be decreasing in the direction of $\widehat{w}_x$, 
moving from $x > 0$ (convex) to $x < 0$ (concave).
Similarly, we find that
$$
\widehat{w}_y \; = \;  -
\begin{bmatrix}
0  \\
m \\
-m_3 \\
0
\end{bmatrix} 
\quad \mbox{and} \quad
\widehat{w}_d \; = \;  
\begin{bmatrix}
0  \\
0 \\
0 \\
-1
\end{bmatrix} 
$$
are outward normals on the boundaries $y=0$ and $d=0$, respectively.

\begin{theorem}
There exists a smooth manifold $\M \subset \Nvex$ with boundary for which the gradient vector field $X(z)$ exists and points outward on the boundary of $\M$.
\label{thm:gradboundary}
\end{theorem}

\pf  
In terms of the variables $x, y,$ and $d$, the gradient of $U$ is
$$
\nabla U(z) \; = \;  
\begin{bmatrix}
-m m_1 r_{12}^{-3} x - m_1 m_3 r_{13}^{-3} (x+y)   \\[0.07in]
-m m_3 r_{23}^{-3} y - m_1 m_3 r_{13}^{-3} (x+y)  \\[0.07in]
-m m_1 r_{12}^{-3} x + m m_3 r_{23}^{-3} y  \\[0.07in]
-m d \left( m_1 r_{12}^{-3}  + m_3 r_{23}^{-3} + m r_{24}^{-3} \right)
\end{bmatrix} .
$$
Pick a point $z_0 = (0, y_0, d_0) \neq P_i$ on the boundary $x=0$ (bodies 1, 2, and 4 collinear).  Then the derivative of $U$ in the direction $\widehat{w}_x$ is
$$
DU(z_0) \widehat{w}_x \; = \;  (\nabla U(z_0))^{\footnotesize{T}} \widehat{w}_x \; = \;  m m_1 m_3 y_0 \left( r_{13}^{-3} - r_{23}^{-3} \right) \; > \; 0 ,
$$
where the inequality follows because $r_{23}$ is the length of the hypotenuse of the right triangle formed by bodies $1, 2,$ and~$3$.  Likewise, for a
point $z_0 = (x_0, 0, d_0) \neq P_i$ on the boundary $y=0$ (bodies 2, 3, and 4 collinear), 
$$
DU(z_0) \widehat{w}_y \; = \;  (\nabla U(z_0))^{\footnotesize{T}} \widehat{w}_y \; = \;  m m_1 m_3 x_0 \left( r_{13}^{-3} - r_{12}^{-3} \right) \; > \; 0 ,
$$
where the inequality follows because $r_{12}$ is the length of the hypotenuse of the right triangle formed by bodies $1, 2,$ and~$3$. 
This shows that the gradient vector field $X|_{\Nvex}$ points outward along the boundaries $x = 0$ and $y = 0$ (excluding the collision points).

We now prove the theorem by contradiction.  Let $\delta > 0$ be a small number.  We claim that it is possible to shrink the boundary of $\Nvex$ slightly
to a smooth curve $\gamma$ that contains $x = 0, y = 0,$ and $d = \delta$, smoothly connects these arcs while avoiding the three collision points $P_i$, 
and for which the gradient vector field points outward.  Suppose this was not the case.  Then, since $\Nvex$ is compact, by taking $\gamma$ closer and closer to the 
boundary of $\Nvex$, we can construct a sequence of points $\{z^k\} \subset \Nvex$ converging to $\overline{z} \in \partial \Nvex$ 
and a corresponding sequence of outward pointing normal vectors
$\{\widehat{w}^k\}$ such that $DU(z^k) \widehat{w}^k \leq 0 \; \forall k$.  We may assume that the normal vectors are all of unit length.

Now, $\overline{z}$ cannot lie on the boundary $x=0$ (excluding $P_2$ and $P_3$) because the limiting outward normal vector would be a positive scalar multiple of $\widehat{w}_x$, and 
$DU( \overline{z} ) \widehat{w}_x > 0$.  Likewise, $\overline{z}$ cannot lie on the boundary $y=0$ (excluding $P_1$ and $P_3$) because the limiting outward normal vector 
would be a positive scalar multiple of $\widehat{w}_y$, contradicting the fact that $DU( \overline{z} ) \widehat{w}_y > 0$.  Finally, 
$\overline{z}$ cannot lie on the boundary $d=0$ (excluding $P_1$ and $P_2$) because the limiting outward normal vector 
would be a positive scalar multiple of $\widehat{w}_d$.  Since $r_{24}^k \rightarrow 0$ while all of the other mutual distances do not converge to 0, we see that
$\dsty{\lim_{k \rightarrow \infty}} DU( z^k ) \widehat{w}^k = +\infty$, which is a clear contradiction.

It follows that $\overline{z}$ must be one of the three corner points $P_i$.  Suppose that $\overline{z} = P_3$, which is the corner point at the North pole where
$x = y = 0$ and $r_{13} = 0$ (binary collision).  Since $P_3$ is a corner, the limiting normal vector does not necessarily exist.  Nevertheless, we will show that
$DU( z^k ) \widehat{w}^k > 0$ for some $k$ by showing that $(\nabla U(z^k))^{\footnotesize{T}} \widehat{w}^k > 0$ for $z^k$ sufficiently close to $P_3$. 
Suppose that $\gamma(t)$ is a smooth arc close to $P_3$ monotonically increasing from the boundary $x=0$ to $y=0$.  
Using equations (\ref{eq:newI}) and~(\ref{eq:normalVec}), we compute
\begin{eqnarray*}
\widehat{w}_1 + \widehat{w}_2   & = &  m^2 m_1 m_3 (x+y) \left( d' \, c - d \, c' \right) \\[0.07in]
& = &  - \frac{m m_1 m_3 (x+y)}{d} \left[  m_1 x' (1 - m_3y(x+y)) - m_3 y' (1 - m_1 x (x+y) ) \right]   
\end{eqnarray*}
for the sum of the first two components of the outward normal vector to $\gamma$. 
Since $x' > 0$ and $y' < 0$ on such a $\gamma$, we see that $\widehat{w}_1 + \widehat{w}_2 < 0$ for $x$ and $y$ sufficiently small.
Then we have
$$
(\nabla U(z))^{\footnotesize{T}} \widehat{w}  \; = \;  - m_1 m_3 r_{13}^{-2} (\widehat{w}_1 + \widehat{w}_2) + R \, ,
$$
where $R$ remains bounded as $z^k \rightarrow P_3$.  Since $r_{13}^k \rightarrow 0$, it follows that 
$(\nabla U(z^k))^{\footnotesize{T}} \widehat{w}^k > 0$ for some $k$, which is a contradiction.

The argument for the other two corner points is similar. Although we do not know the limiting normal vector,
we know the signs of the key components of $\widehat{w}^k$ for $k$ sufficiently large.   
For example, if $\overline{z} = P_2 = (0,y_m,0)$, then $x$ and $d$ can be taken sufficiently close to 0 to insure that $\widehat{w}_1 + \widehat{w}_3 < 0$ and $\widehat{w}_4 < 0$.  
Then we find
$$
(\nabla U(z))^{\footnotesize{T}} \widehat{w}  \; = \;  - m m_1 r_{12}^{-3} x (\widehat{w}_1 + \widehat{w}_3) - md(m_1 r_{12}^{-3} + mr_{24}^{-3}) \widehat{w}_4 +  R \, ,
$$
where $R$ remains bounded as $z^k \rightarrow P_2$.  Recall that $r_{12} = \sqrt{x^2 + d^2}$ and $r_{24} = 2d$.
Since each of these mutual distances approaches 0,  it follows that 
$(\nabla U(z^k))^{\footnotesize{T}} \widehat{w}^k > 0$ for some $k$, which is a contradiction.
The argument for $\overline{z} = P_1$ is similar.  This completes the proof.
\enpf

\begin{remark}
We note that $\M$ can be chosen so that there are no critical points of $V$ in $\Nvex - \M$; otherwise we could construct a sequence of critical points converging to 
a point on the boundary of $\Nvex$ corresponding to collision.  But this is impossible by the same argument used in the proof of the theorem.  This also
follows by Shub's Lemma~\cite{shub}.
\end{remark}

\begin{corr}[Existence]  For any choice of positive masses $(m_1, m_3, m_2 = m_4)$, there exists a convex kite central configuration with bodies 2 and~4 equidistant
from the axis of symmetry.
\label{Cor:Exists}
\end{corr}

\pf
Since $\M$ is a compact manifold and $V$ is continuous on $\M$, $V$ must attain a minimum~$z$.  By Theorem~\ref{thm:gradboundary},
the gradient flow points outward on the boundary of $\M$, so $z$ must lie in the interior of $\M$ and is therefore a convex configuration.
Since $z$ is a critical point of $V$, it is also a kite central configuration.  \enpf

\section{Uniqueness of Convex Kites}
\label{Sec:Unique}

In this section we prove that for any choice of positive masses $(m_1, m_3, m_2 = m_4 = m/2)$, there exists a {\em unique} convex kite central configuration with
our prescribed ordering.  Our proof uses the Poincar\'{e}-Hopf Index Theorem, Theorem~\ref{thm:gradboundary}, and formula~(\ref{eq:prodEvals}).

First we review some important concepts from differential topology~\cite{guilpoll, milnor}.  Suppose we have a smooth vector field $\nu: U \mapsto \mathbb{R}^n$ on some open
set $U \subset \mathbb{R}^n$ containing an isolated zero $z$.  The {\em index} of $z$, denoted $i(z)$, is the degree of the map $s \mapsto \nu(s)/||\nu(s)||$, which maps a small ball about $z$
onto the unit sphere.  Since $z$ is isolated, the domain of this map can be chosen sufficiently small so as not to contain any other zeros than $z$.
For vector fields in the plane, $i(z)$ is obtained by counting the integer number of rotations of $\nu(s)$ while traversing a
simple closed loop about $z$ in the counterclockwise direction.  Rotations of the vector field in the counterclockwise direction add $+1$ to the index, while rotations
in the clockwise direction contribute $-1$.

The index $i(z)$ for a vector field on a manifold can be obtained by using local parametrizations.  Suppose now that $\nu$ is a smooth vector field on some manifold $\M$ with isolated
zero $z$ and that $\phi: U \mapsto \M$ is a local parametrization (diffeomorphism) of a neighborhood of $z$.  Then we define $i(z)$ to be the index of 
the corresponding vector field $d\phi^{-1} \circ \nu \circ \phi$ (the pullback) on $U$ at the isolated zero $\phi^{-1}(z)$.  For our purposes, the salient fact is that
the index of a nondegenerate critical point is $+1$ if the determinant of the linear map $d\nu(z): T_z \M \mapsto T_z \M$ is positive (see Lemma~5 in Section~6 of~\cite{milnor}).
Any nondegenerate source or sink has an index of $+1$.

The following classical theorem, adapted to manifolds with boundary~\cite{milnor}, is the key tool in our proof.

\begin{theorem}[Poincar\'{e}-Hopf Index Theorem]
Let ${\cal{M}}$ be a smooth compact manifold with boundary and let $X$ be a vector field on ${\cal{M}}$ with isolated zeros
such that $X$ points outwards along the entire boundary.  Then the sum of the indices of the zeros is equal to the Euler characteristic of~${\cal{M}}$.
\end{theorem}

We will apply this theorem to the smooth manifold $\M$ given by Theorem~\ref{thm:gradboundary} 
and the vector field $X(z)$ from Lemma~\ref{lemma:gradfield}.  First we show that the index of {\em any} zero of $X$ has index $+1$.

\subsection{Reducing the calculation to two variables}
\label{subsec:reductwovar}

By Lemma~\ref{lemma:invBperp}, the modified Hessian $M^{-1}H(z)$ and $dX(z)$ are the same when restricted to $T_z \N$.  Thus, the 
determinant of $dX(z)$, when viewed as a linear map on $T_z \N$, is equal to the product of the nontrivial eigenvalues of $M^{-1}H(z)$.
Recall from formula~(\ref{eq:prodEvals}) that this product is given by
\begin{equation}
\frac{1}{2} \left[  ( \, \mbox{tr} (A) \, )^2 - \mbox{tr}( A^2 ) \right]  - \lambda \mbox{tr}(A) + 3\lambda^2 \, ,
\label{eq:indexA}
\end{equation}
where $A = M^{-1}D^2U(z)$.
Due to the homogeneity of the potential function $U$, we can scale the variables by a common factor without changing the signs of the eigenvalues.
While this scaling moves $z$ off the manifold $\M$ (since $I$ no longer equals 1), it does not alter the index of~$z$.

Specifically, suppose we replace $z = (a,b,c,d)$ by $kz = (ka, kb, kc, kd)$, for some positive factor $k$.  Then all of the mutual distances
$r_{ij}$ scale by a factor of $k$ and according to equation~(\ref{eq:lambda}), $\lambda$ scales by $k^{-3}$.  If we set $\wl = k^{-3} \lambda$, 
then the new modified Hessian is
$$
M^{-1} H(kz) \; = \; M^{-1} D^2U(kz)+ \wl \; = \;  k^{-3} \left( M^{-1} D^2U(z) + \lambda \right) \; = \;  k^{-3} M^{-1} D^2H(z)  \, .
$$
It follows that the eigenvalues are scaled by the positive factor $k^{-3}$.  Moreover, the sign of~(\ref{eq:indexA}), and
therefore the index of $z$, is also unaffected by the scaling, since $\mbox{tr}(A)$ scales by the positive factor $k^{-3}$ and
$\mbox{tr}(A^2)$ scales by $k^{-6}$.

Next we observe that all of the entries in $A$ (including the mutual distances) are functions of $x = a+c$, $y = b-c$, $x+y = a+b$, and $d$.  
Thus, the nontrivial eigenvalues are determined by the values of $(x,y,d)$.   The variable $d$ (assumed positive) can be eliminated by 
scaling all variables by the factor $k = 1/d$:
$$
(x, \, y, \, d) \mapsto (\wx, \, \wy, \, 1) \, ,
$$ 
where $\wx = x/d$ and $\wy = y/d$.
As discussed above, this change of variables does not alter the signs of the eigenvalues.  In essence,
we have shifted the configuration $z$ to the right by $c$ units and rescaled it by the factor $k = 1/d$.  In $(\wx, \wy)$ coordinates, 
the positions of the bodies are $q_1 = (\wx,0), q_2 = (0,1), q_3 = (-\wy,0),$ 
and $q_4 = (0,-1)$, with the center of mass at $(\wc,0)$, where $\wc = m_1 \wx - m_3 \wy$.
Although this moves $z$ off the normalized configuration space, the signs of the eigenvalues (and hence the index)
are unchanged.  The configuration is convex when $\wx \wy > 0$ and concave when $\wx \wy < 0$.

Letting $\wa = a/d$ and $\wb = b/d$, we have
\begin{eqnarray*}
\wa \; = \; \wx - \wc \; = \;  \wx - m_1 \wx + m_3 \wy & = &  m_3(\wx + \wy) + m \wx \; \mbox{ and}\\
\wb \;  = \; \wy + \wc \;  = \;  \wy + m_1 \wx - m_3 \wy & = &  m_1(\wx + \wy) + m \wy \, .
\end{eqnarray*}
Using these, the central configuration equations (\ref{eq:cc1}), (\ref{eq:cc2}), and~(\ref{eq:lambda}) become
\begin{eqnarray}
m \wx(r_{12}^{-3} - \wl)  & = &  m_3 (\wx + \wy)(\wl - r_{13}^{-3}) \label{eq:cc1-new} \\[0.07in]
m \wy(r_{23}^{-3} - \wl)  & = &  m_1 (\wx + \wy)(\wl - r_{13}^{-3}) \label{eq:cc2-new} \\[0.07in]
\wl & =  &  m_1 r_{12}^{-3} + m_3 r_{23}^{-3} + m r_{24}^{-3}  \label{eq:lambda-new} \, ,
\end{eqnarray}
where $\wl = d^3 \lambda$ and the mutual distances are now
\begin{equation}
r_{12} = \sqrt{\wx^2 + 1}, \quad  r_{23} = \sqrt{\wy^2 + 1},  \quad  r_{13} = \wx + \wy, \quad \mbox{and } r_{24} = 2.
\label{eq:mds-new}
\end{equation}
As a check on our work, we note that equations (\ref{eq:cc1-new}), (\ref{eq:cc2-new}), and~(\ref{eq:lambda-new}) can also be derived by
substituting $q_1 = (\wx,0), q_2 = (0,1), q_3 = (-\wy,0),$ and $q_4 = (0,-1)$ into equation~(\ref{cc:maineq}).

In $(\wx,\wy)$ coordinates, the key matrix $A = M^{-1} D^2U(z)$ becomes
\begin{equation}
\begin{bmatrix}
\zeta_1  & 2m_3 r_{13}^{-3} &  m r_{12}^{-5} \left( 2\wx^2 - 1 \right) &  3m \wx r_{12}^{-5}    \\[0.07in]
2m_1 r_{13}^{-3} &  \zeta_2 & -m r_{23}^{-5} \left( 2\wy^2 - 1 \right)  &  3m \wy r_{23}^{-5}    \\[0.07in]
m_1 r_{12}^{-5} \left( 2\wx^2 - 1 \right)   & -m_3 r_{23}^{-5} \left( 2\wy^2 - 1 \right)  & \zeta_3  &  3m_1 \wx r_{12}^{-5} - 3m_3 \wy r_{23}^{-5}   \\[0.07in]
3m_1 \wx r_{12}^{-5} &  3m_3 \wy r_{23}^{-5} &  3m_1 \wx r_{12}^{-5} - 3m_3 \wy r_{23}^{-5} & \zeta_4   
\end{bmatrix}
\, ,
\label{matrix:redHess}
\end{equation}
with
\begin{eqnarray}
\zeta_1 = m r_{12}^{-5} \left( 2 \wx^2 - 1 \right) + 2m_3r_{13}^{-3},  & \quad &  \zeta_2 = m r_{23}^{-5} \left( 2 \wy^2 - 1 \right) + 2m_1r_{13}^{-3},  \label{eq:trace1-new} \\[0.07in] 
\zeta_3 = m_1 r_{12}^{-5} \left( 2 \wx^2 - 1 \right) + m_3 r_{23}^{-5} \left( 2 \wy^2 - 1 \right),   & \quad & \zeta_4 = 3m_1r_{12}^{-5} + 3m_3r_{23}^{-5} + 3mr_{24}^{-3} - \wl \, .
\label{eq:trace2-new}
\end{eqnarray}

\subsection{The reduced configuration space $\C$}
\label{subsec:reducedC}

Note that equation~(\ref{eq:lambda-dist}) still holds in our new coordinates, with $\lambda$ replaced by $\wl$ and the mutual distances 
given by equation~(\ref{eq:mds-new}).  Using this expression for $\wl$,
we can derive formulas for the mass ratios that only depend on $\wx$ and~$\wy$.
Equations (\ref{eq:cc1-new}) and~(\ref{eq:cc2-new}) imply
\begin{eqnarray}
\frac{m_1}{m} & = &   \frac{\wy}{\wx + \wy} \cdot  \frac{r_{23}^{-3} - \wl}{\wl - r_{13}^{-3}} \; = \; 
                                   \frac{\wy}{\wx + \wy} \cdot  \frac{r_{23}^{-3} - r_{24}^{-3}}{r_{12}^{-3} - r_{13}^{-3}} \label{eq:massratio1} \\[0.07in]
\frac{m_3}{m} & = &   \frac{\wx}{\wx + \wy} \cdot  \frac{r_{12}^{-3} - \wl}{\wl - r_{13}^{-3}} \; = \; 
				  \frac{\wx}{\wx + \wy} \cdot  \frac{r_{12}^{-3} - r_{24}^{-3}}{r_{23}^{-3} - r_{13}^{-3}} \; . \label{eq:massratio2}
\end{eqnarray}

Without loss of generality, we can make one final reduction by assuming that $r_{12} \geq r_{23}$, which is equivalent to $|\wx| \geq |\wy|$.
This implies that the longest sides of a convex kite are between bodies 1 and 2, and bodies 1 and 4.  If this was not the case,
then we could rotate the configuration by $180^{\circ}$ and relabel the bodies.  Alternatively, we note that
equations (\ref{eq:cc1-new}) and~(\ref{eq:cc2-new}) possess a symmetry upon interchanging the variables $\wx$ and $\wy$ while also interchanging
the parameters $m_1$ and $m_3$.
Thus, if $(\wx, \wy)$ is a kite central configuration for the masses $(m_1, m_3, m)$, then so is $(\wy, \wx)$, with masses $(m_3, m_1, m)$.  The index of each 
of these critical points is identical since the potential function is invariant under this symmetry.

\begin{figure}[t]

\begin{center}
\includegraphics[height = 320bp]{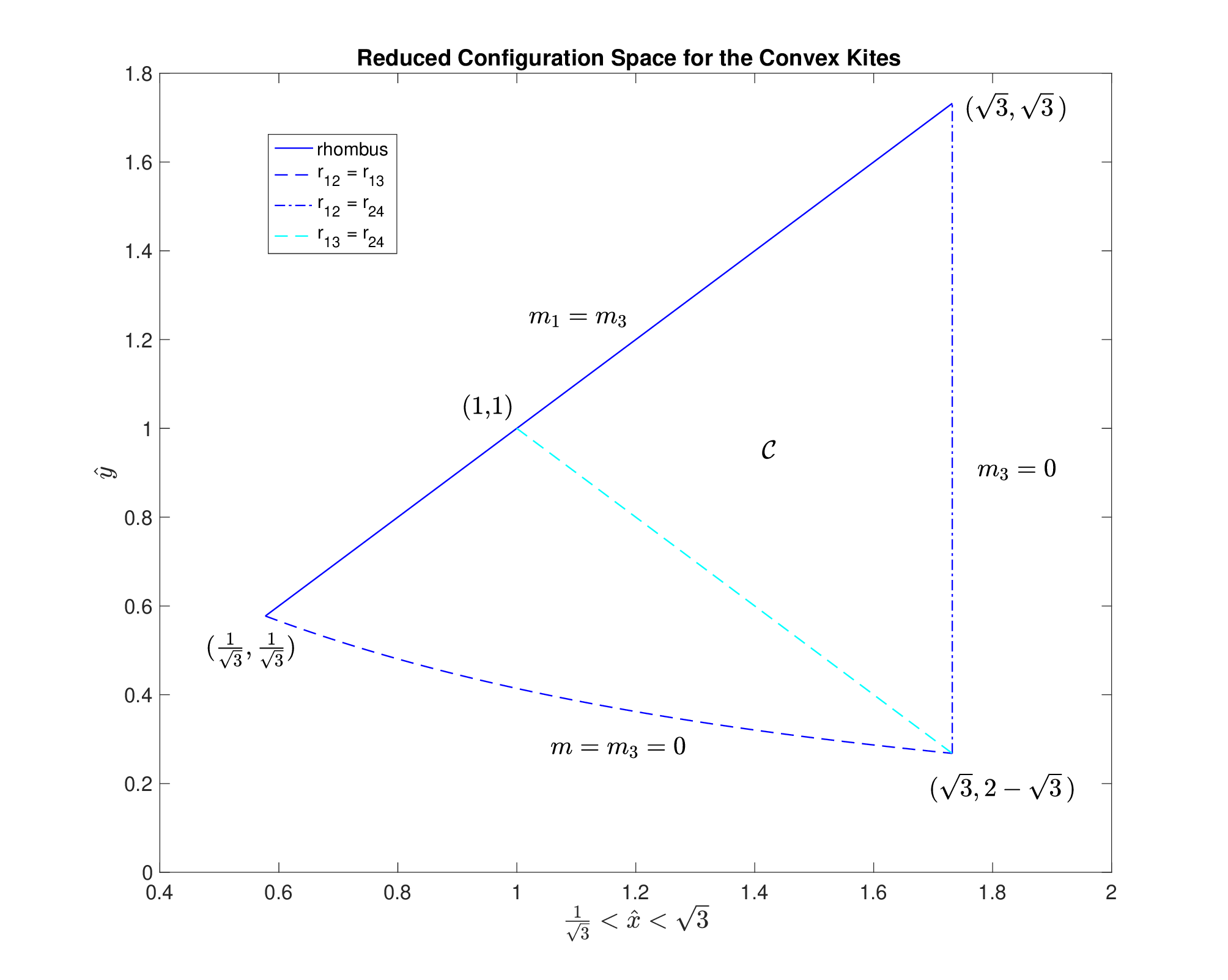}
\end{center}
\vspace{-0.25in}
\caption{The reduced configuration space $\C$ for the convex kites, defined by $r_{23} \; \leq \; r_{12}  \; <  \; r_{13}, r_{24}$.
Configurations on the lower and right boundaries, respectively, correspond to c.c.'s of the $1+3$-body and restricted four-body problems, respectively.}
\label{Fig:ReducedDomain}
\end{figure}

Focusing on the convex case, we may assume that $\wx \geq \wy > 0$.  Then, 
positivity of the masses and equations (\ref{eq:massratio1}) and~(\ref{eq:massratio2}) imply the following relations among the mutual distances: 
$$
r_{23} \; \leq \; r_{12}  \; <  \; (\wl)^{-1/3} \;  < \; r_{13}, r_{24} \, .
$$
The diagonals of the kite must be longer than all of the exterior sides, a fact that holds for any convex central configuration, not just kites~\cite{schmidt}.
The mutual distance inequalities describe a bounded subset $\C$ of the $\wx \wy$-plane (see Figure~\ref{Fig:ReducedDomain}).
This set has three boundary curves, the first of which is in $\C$:
\begin{itemize}
\item[{\bf (i)}]  $r_{12} = r_{23} \; \Longrightarrow \;  \wy = \wx, \, 1/\sqrt{3} < \wx < \sqrt{3} \mbox{ and }  m_1 = m_3$ (rhombus configuration),

\item[{\bf (ii)}]  $r_{13} = r_{12} \; \Longrightarrow \;  \wy = -\wx + \sqrt{\wx^2 + 1} \mbox{ and }  m = m_3 = 0$ (bodies 2, 3, and 4 equidistant from body 1),

\item[{\bf (iii)}]  $r_{24} = r_{12} \; \Longrightarrow \;  \wx = \sqrt{3} \mbox{ and }  m_3 = 0$ (bodies 1, 2, and 4 form an equilateral triangle).
\end{itemize}
Configurations on boundary~{\bf (ii)} are c.c.'s of the $1+3$-body problem, with one large mass $m_1$ and three infinitesimal masses.
Configurations on boundary~{\bf (iii)} are c.c.'s of the restricted four-body problem, with three ``primaries'' located at the vertices of an equilateral triangle
and a fourth infinitesimal mass on the axis of symmetry of the kite.

We will call $\C$ the {\em reduced configuration space}, defined as
$$
\C  \; = \; \{ (\wx, \wy) \in \mathbb{R}^2: \frac{1}{\sqrt{3}} < \wx < \sqrt{3}, \, -\wx + \sqrt{\wx^2 + 1} \, < \, \wy \, \leq \, \wx \}  \, .
$$
For each $(\wx, \wy)$ in $\C$, the four bodies with positions $q_1 = (\wx,0), q_2 = (0,1), q_3 = (-\wy,0), q_4 = (0,-1)$ form a
convex kite central configuration with positive masses.  The values of the masses are uniquely determined by
equations (\ref{eq:massratio1}), (\ref{eq:massratio2}), and the normalization $m_1 + m_3 + m = 1$.
Any point in the first quadrant that is not contained in $\C$ corresponds to a convex kite c.c.~with 
a negative (or zero) mass, or a configuration satisfying $r_{23} > r_{12}$.

\subsection{The index of a convex kite central configuration}

\begin{lemma}
Suppose that $z \in \Nvex$ is a convex kite central configuration with positive masses.  Then, when viewed as a rest point of the gradient vector field $X$, 
the index of $z$ is $+1$.
\label{lemma:index}
\end{lemma}

\pf
Recall that the index of a rest point $z$ is $+1$ if the product of the nontrivial eigenvalues of $M^{-1} H(z)$ is positive.  
By the arguments in Sections \ref{subsec:reductwovar} and~\ref{subsec:reducedC}, it suffices to show that 
\begin{equation}
\frac{1}{2} \left[  ( \, \mbox{tr} (A) \, )^2 - \mbox{tr}( A^2 ) \right]  - \wl \mbox{tr}(A) + 3\wl^2 \; > \; 0
\label{indexCond}
\end{equation}
for any point $(\wx, \wy) \in \C$, where $A = M^{-1}D^2U(z)$ is given by~(\ref{matrix:redHess}) and $\wl$ is determined by equation~(\ref{eq:lambda-new}).

To simplify our notation, we will drop the ``$\; \hat{\;} \;$'' from $\wx$ and $\wy$ for the remainder of the proof. 
Let
$$
\alpha_1 = 2x^2 - 1, \quad  \alpha_2 = 2y^2 - 1, \quad  \beta_1 = 4x^2 + 1, \quad \beta_2 = 4y^2+ 1 \, .
$$
After a lengthy computation, condition~(\ref{indexCond}) is equivalent to showing that 
\begin{eqnarray}
& & m^2 \left( \casefrac{1}{64} + \casefrac{1}{8} r_{12}^{-5} \alpha_1 + \casefrac{1}{8} r_{23}^{-5} \alpha_2 + r_{12}^{-5} r_{23}^{-5} \alpha_1 \alpha_2 \right) \nonumber
+  m m_1 \left( \casefrac{1}{8} r_{12}^{-3} + \casefrac{3}{4} r_{12}^{-5} x^2 + 2 r_{12}^{-5} r_{13}^{-3} \alpha_1 - r_{12}^{-8} \beta_1 + \casefrac{1}{4} r_{13}^{-3}  \right)  \\[0.07in]
& &  + \, m m_3 \left( \casefrac{1}{8} r_{23}^{-3} + \casefrac{3}{4} r_{23}^{-5} y^2 + 2 r_{23}^{-5} r_{13}^{-3} \alpha_2 - r_{23}^{-8} \beta_2 + \casefrac{1}{4} r_{13}^{-3}  \right) 
+  9 m_1m_3  r_{12}^{-5} r_{23}^{-5} (x+y)^2  \; > \; 0 \, .  \label{eq:condwMass}
\end{eqnarray}
If we divide this inequality by $m^2$, we can eliminate the masses entirely by substituting in 
formulas (\ref{eq:massratio1}) and (\ref{eq:massratio2}).  This is an important step because the left-hand side of
inequality~(\ref{eq:condwMass}) vanishes on the lower boundary of $\C$, where $m = m_3 = 0$.  Performing this substitution and then
multiplying through by the denominator $64 r_{12}^5 r_{23}^5 r_{13} (r_{13}^3 - r_{12}^3) (r_{13}^3 - r_{23}^3)$ (which is strictly positive on~$\C$)
leads to the inequality $F(x,y) > 0$, where $F$ is the sum of the four functions
\begin{eqnarray*}
F_1  &  =  &  r_{13}(r_{13}^3 - r_{12}^3)(r_{13}^3 - r_{23}^3)(r_{12}^5 + 8\alpha_1)(r_{23}^5 + 8\alpha_2),  \\[0.07in]
F_2  &  =  &  y \, r_{23}^2(r_{13}^3 - r_{23}^3)(8 - r_{23}^3) [r_{13}^3 (r_{12}^5 + 6x^2 r_{12}^3  - 8\beta_1) + 2r_{12}^3(r_{12}^5 + 8\alpha_1) \, ],   \\[0.07in]
F_3  &  =  &  x \, r_{12}^2(r_{13}^3 - r_{12}^3)(8 - r_{12}^3) [r_{13}^3 (r_{23}^5 + 6y^2 r_{23}^3 - 8\beta_2) + 2r_{23}^3(r_{23}^5 + 8\alpha_2) \, ],  \\[0.07in] 
F_4  &  =  &  9 xy \, r_{13}^7 (8 - r_{12}^3)(8 - r_{23}^3)  \, .
\end{eqnarray*}
We will use Gr\"{o}bner bases and symmetry to show that $F(x,y) > 0$ for all $(x,y) \in \C$.  All calculations are symbolic only (no numerical
estimates required) and were performed with SageMath~\cite{sage}.



Recall that $r_{12} = \sqrt{x^2+1}, r_{23} = \sqrt{y^2+1}$, and $r_{13} = x+y$.  As expected, the function $F$ is symmetric, that is, $F(x,y) = F(y,x)$.  
The first step is to bound $F$ below by a symmetric polynomial $P(x,y)$ using quadratic bounds for $r_{12}$ and $r_{23}$.
It is straight-forward to check that
\begin{eqnarray}
\casefrac{1}{6} x^2 + \casefrac{1}{3} x + \casefrac{8}{9} \quad \leq & r_{12} &  \leq \quad  \casefrac{1}{6} x^2 + \casefrac{1}{3} x + \casefrac{17}{18}, \label{quadest1} \\[0.07in]
\casefrac{1}{6} y^2 + \casefrac{1}{3} y + \casefrac{8}{9} \quad \leq & r_{23} &  \leq \quad  \casefrac{1}{6} y^2 + \casefrac{1}{3} y + \casefrac{17}{18}  \label{quadest2}
\end{eqnarray}
for $(x,y) \in \C$. To obtain the polynomial $P$, we expand each $F_i$ function, simplifying any even powers of $r_{12}$ and $r_{23}$, and writing
all odd powers as $r_{12} \cdot (r_{12})^{2j}$ and $r_{23} \cdot (r_{23})^{2j}$ before simplifying. After 
carefully grouping those terms with $r_{12}, r_{23},$ and $r_{12} r_{23}$ by sign,
we then apply inequalities (\ref{quadest1}) and~(\ref{quadest2}) to bound each $F_i$ below by polynomials $P_i(x,y)$.  
This produces a polynomial $P(x,y) = P_1 + \cdots + P_4$ with rational coefficients that is a lower bound for $F$ on $\C$.
$P$ is a degree 19 symmetric polynomial with 165 terms.

Next we show that $P$ has no critical points in the interior of $\C$.  We take advantage of the symmetry by introducing the variables
$$
\tau \; = \;  x + y \quad \mbox{ and } \quad \rho = xy \, .
$$
Using Gr\"{o}bner bases, we convert $P$ into a polynomial in the variables $(\tau,\rho)$:
$$
\begin{array}{c}
SP(\tau,\rho) =   3\tau^{13} + 24 \tau^{12} - 758\tau^{11}\rho^2 + 370\tau^{11}\rho + 87\tau^{11} + 36\tau^{10}\rho^3 + \cdots +  \\[0.07in]
    41728\tau \rho - 11812 \tau - 3456\rho^7 + 17280\rho^6 - 34560\rho^5 + 34560\rho^4 - 17280\rho^3 + 3456\rho^2 .
\end{array}
$$
$SP$ is a degree 13 polynomial with integer coefficients and 79 terms.  The specific relationship between $P$ and $SP$ is given by
$$
SP(\tau(x,y), \rho(x,y)) \; = \;  -108P(x,y) \, .
$$
Solving the system $\{\partial SP/\partial \tau = 0, \partial SP/\partial \rho =0 \}$, we use Gr\"{o}bner bases and a lex order to eliminate $\rho$.  After dividing out the
trivial root at $\tau = 0$, we obtain a single variable polynomial in $\tau$ of degree 112.  Using the command \verb+number_of_roots_in_interval+ in Sage
(which executes a Sturm algorithm), we confirm that this polynomial has no roots in the interval $2/\sqrt{3} < \tau <  2\sqrt{3}$.  By the chain rule,
it follows that $P$ has no critical points on the interior of~$\C$.

Therefore, $P(x,y)$ must attain it's minimum value on the boundary of $\C$.  It is straight-forward to check that $P > 0$ on each boundary since $P$ reduces to
a function of one variable in each case.  On the diagonal $y = x$, $P$ is strictly increasing with a minimum value of approximately $162.38$ at $x = 1/\sqrt{3}$.
On the vertical boundary $x=\sqrt{3}$, $P$ is also strictly increasing with a minimum of approximately $12,\!818$ at $y = 2 - \sqrt{3}$.
Finally, for the lower boundary, in order to handle the $\sqrt{x^2 + 1}$ term, we introduce the auxiliary variable $u$ and consider the ideal in $\mathbb{Q}[u,x]$ generated by the 
polynomials $P(x,y=-x+u)$ and $u^2 - x^2 - 1$.  Computing a Gr\"{o}bner basis for this ideal to eliminate $u$, we obtain a polynomial in $x$ that has
no roots for $1/\sqrt{3} < x < \sqrt{3 \,}$.  Since $P(1,-1+\sqrt{2}) > 0$, we conclude that $P > 0$ on the entire lower boundary. (The minimum value is the same as the
value attained on the diagonal.)  This shows that $P$ is strictly positive on $\C$.

Although the specific location of the critical point $z$ is unknown, we do know that it must correspond to some point $(\wx, \wy)$ in the reduced configuration space~$\C$.
Since $F(\wx, \wy) > P(\wx, \wy) > 0$, it follows that the index of $z$ is $+1$.
\enpf

\begin{remark}
Using MATLAB~\cite{matlab}, we evaluated $F$ on a $300 \times 300$ grid on $\overline{\C}$ to check numerically that $F > 0$.
The minimum value of $F$ is approximately $342.71$, attained at the point $(1/\sqrt{3}, 1/\sqrt{3})$, which is the same location
for the minimum of $P$.  It appears that $F$ is an increasing function of $y$ on $\overline{\C}$.
\end{remark}

\subsection{Proof of Main Theorem}

\begin{theorem}
Given four positive masses $m_1, m_3, m_2 = m_4$ and any ordering with bodies 2 and 4 positioned opposite each other, there exists a unique convex kite central configuration.
\label{thm:main}
\end{theorem}

\pf
Without loss of generality, we can assume that the bodies are ordered consecutively in the counterclockwise direction.  A convex kite central configuration with this ordering corresponds
to a critical point $z \in \Nvex$ that is a rest point (zero) of the gradient vector field $X$ defined in Lemma~\ref{lemma:gradfield}.   
Let $Z \subset \Nvex$ be the set of all zeros of $X$ contained in the space of normalized convex configurations.
By Corollary~\ref{Cor:Exists}, this set is nonempty.

By Theorem~\ref{thm:gradboundary}, there exists a smooth manifold $\M \subset \Nvex$ with boundary for which $X$ points outward on the boundary.
We have $Z \subset M$ by construction.  Since
$\M$ is diffeomorphic to a closed disk, it has Euler characteristic~1.
Each $z \in Z$ is an isolated zero of the vector field $X$ with index equal to $+1$.
The zeros are isolated because there are a finite number of them~\cite{hamp-rick}.   Let $N$ be the number of
elements in $Z$, that is, the total number of equivalence classes of convex kite central configurations with our prescribed ordering.  
Then, by the Poincar\'{e}-Hopf Index Theorem, we have
$$
1 \; = \; \chi( \M )  \; = \;  \sum_{z \in Z}^N  \, i(z)  \; = \;  \sum_{z \in Z}^N  \, 1   \, ,
$$
from which it follows quickly that $N = 1$.
\enpf

\begin{corr}
A convex kite central configuration $z \in \N$ is a nondegenerate local minimum of $V = U|_{\N}$.
\end{corr}

\pf
Without loss of generality, we can assume that $z \in M$. Nondegeneracy follows from Lemma~\ref{lemma:index}.  
If $z$ were not a local minimum, then there would exist an additional critical point $z^\ast \in M$ corresponding to
the minimum value of $V$ on $\M$.  But then $z^\ast$ would be a convex kite central configuration with the same ordering as $z$, contradicting uniqueness.
\enpf

\begin{remark} 
If we expand the configuration space to include {\em all} convex configurations (not just kites), then the
convex kite c.c.'s are local minima of the full four-body potential function $U(q)$ 
restricted to the mass ellipsoid $I(q) = 1$.  This follows from the same argument used in the proof of the corollary along with the fact that
any convex central configuration with masses $m_1, m_3, m_2 = m_4$ must necessarily be a kite~\cite{albouy}. 
\end{remark}

\begin{remark}
One advantage of using the Poincar\'{e}-Hopf Index Theorem is that computing the Poincar\'{e} index is easier than computing the Morse index or proving
that $z$ is a minimum of the restricted potential function.  This should prove fruitful for problems set in higher dimensional spaces.
\end{remark}

\section{Concave Kites}

In this section we discuss the concave kite central configurations.  Much of our work is numerical in nature and concurs with the results of Leandro~\cite{Leandro}.
The setup is similar to the convex setting except that now $y = b - c$ is negative.
Positivity of the masses leads to two different types of concave configurations.  We find a degenerate family of concave kites and 
a mass map that is two-to-one.

Without loss of generality, we assume that the third body ($m_3$) lies in the interior of the triangle formed by the other three
bodies and that bodies 1, 2, and~4 are arranged in counterclockwise order (see right plot in Figure~\ref{Fig:setup}).  
Recall that $x = a + c$, $y = b - c$, $r_{13} = x+y$, and $r_{24} = 2d$.
Central configurations of this type lie in the interior of
$$
\Ncav \; = \;  \{ z \in \N: y \leq 0, \, x + y \geq 0, \mbox{ and } d \geq 0 \} .
$$
Unlike the convex case, the gradient flow points {\em inwards} on the boundary $y = 0$ of $\Ncav$, so the Poincar\'{e}-Hopf Index Theorem
does not apply.  In this case, the minimum value of $U|_{\Ncav}$ may be located on the boundary $y = 0$.
In contrast to the convex setting,  there are values of the masses for which there does not exist a concave kite central configuration (although
$m_2 = m_4$, the central configuration need not be symmetric).  This fact was previously noted in~\cite{albouy}.

\subsection{The reduced configuration space $\Cv_1 \cup \Cv_2$}

The reduction and change of variables from $z = (a,b,c,d)$ to $(\wx, \wy)$, as well as the formulas for the mass ratios given in (\ref{eq:massratio1}) and~(\ref{eq:massratio2}),
are still valid in the concave setting.  However, since $\wy < 0$, we obtain different relations between the mutual distances.  Positivity of the masses and equation~(\ref{eq:Dziobek})
imply that 
$$
r_{23}, \, r_{13} \; < \;  \lambda^{-1/3} \; < \;  r_{12}, \, r_{24} \, .
$$
This means that each edge of the exterior triangle is larger than any of the interior edges.  In addition, equation~(\ref{eq:massratio2}) implies that
$r_{12} > r_{24}$ if and only if $r_{23} > r_{13}$.  In other words, each of the congruent legs of the exterior triangle is larger than the base of the triangle 
if and only if each of the congruent interior edges is larger than the third interior edge.

\begin{figure}[t]

\begin{center}
\includegraphics[height = 350bp]{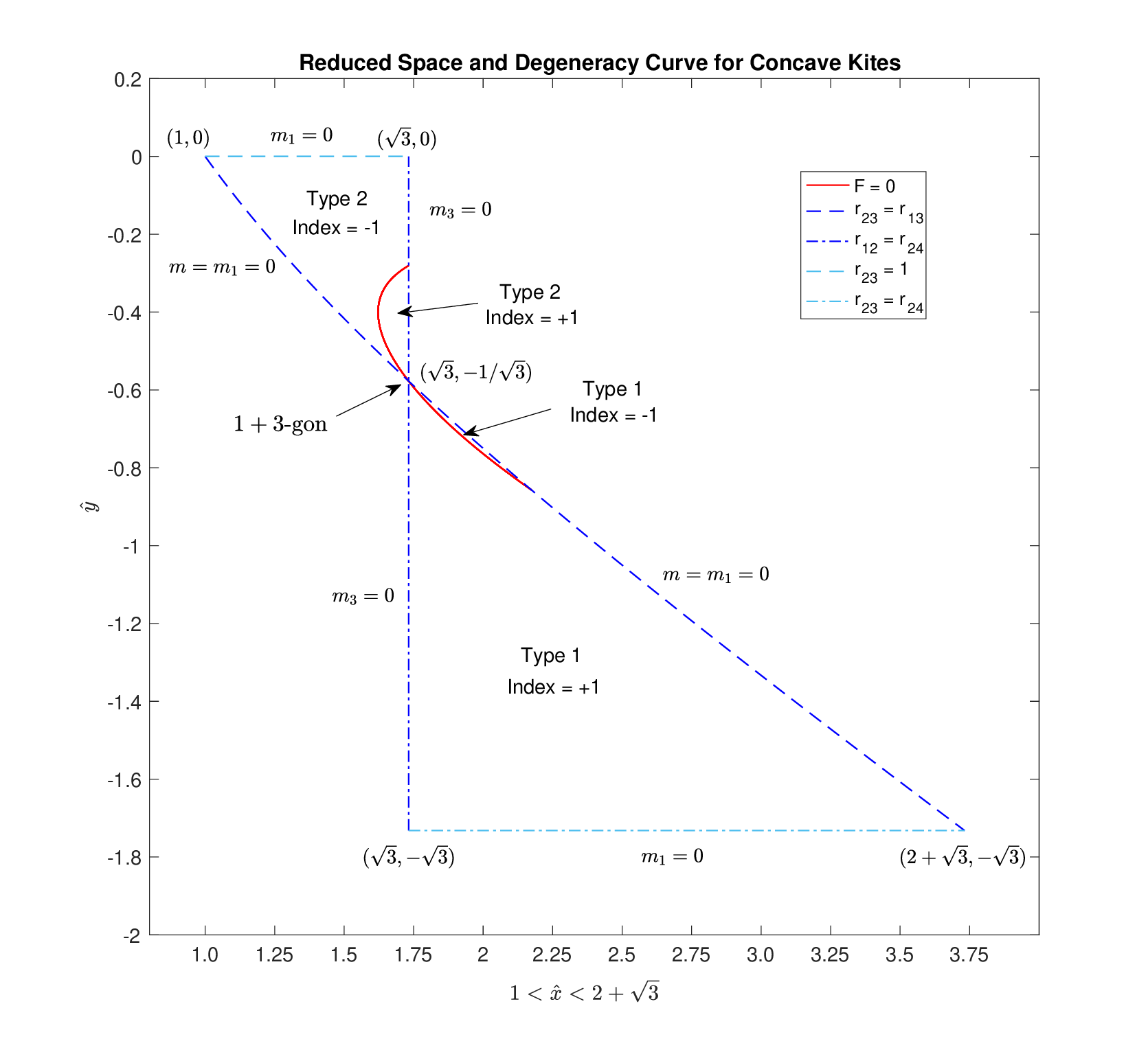}
\end{center}
\vspace{-0.35in}
\caption{The reduced configuration space for the different types of concave kite central configurations.  The red curve consists of degenerate central configurations.}
\label{Fig:ConcaveDomain}
\end{figure}

Thus we have two types of concave kite configurations.  We will refer to the concave configurations
satisfying $r_{12} > r_{24}$ (and consequently $r_{23} > r_{13}$) as {\em Type 1}.   Using~(\ref{eq:mds-new}) and the defining inequalities, Type~1 configurations lie in
the space
$$
\Cv_1  \; = \; \{ (\wx, \wy) \in \mathbb{R}^2:  \sqrt{3} < \wx < 2 + \sqrt{3} \, , \, -\sqrt{3} \, < \, \wy \, < \, -\frac{1}{2} \left( \wx - \frac{1}{\wx} \right) \}  \, .
$$
Concave configurations with $r_{12} < r_{24}$ (and consequently $r_{23} < r_{13}$) will be denoted as {\em Type 2}.   These configurations lie in the space
$$
\Cv_2  \; = \; \{ (\wx, \wy) \in \mathbb{R}^2:  1 < \wx <  \sqrt{3} \, , \,   -\frac{1}{2} \left( \wx - \frac{1}{\wx} \right)  \, < \, \wy \, < \, 0   \}  \, 
$$
(see Figure~\ref{Fig:ConcaveDomain}). 
The closures of these two sets intersect in a special point at $(\sqrt{3}, -1/\sqrt{3})$, which corresponds to the well-studied $1+3$-gon 
central configuration.  This configuration satisfies $r_{12} = r_{24} = 2$, $r_{23} = r_{13} = 2/\sqrt{3}$, and $m_1 = m/2 = m_2 = m_4$.
It consists of three equal mass at the vertices of an equilateral triangle and a body of arbitrary mass located at the center.   This is the only kite configuration
with an arbitrary mass;  each configuration in $\C$, $\Cv_1$, or $\Cv_2$ has a unique value of the normalized masses that makes it central.  
We note the similarity of our bow-tie shaped configuration space in Figure~\ref{Fig:ConcaveDomain} 
with the sets $\Pi_1$ and $\Pi_3$ in Fig.~4 of~\cite{Leandro}.

As with the convex case, the boundary curves (dashed curves in Figure~\ref{Fig:ConcaveDomain}) correspond to central configurations where one or more masses vanish.  These
configurations are central configurations of either the restricted four-body problem or the $1+3$-body problem.  For example, 
the boundary $\wy = 0$ contains solutions to the restricted four-body problem, but with the three primaries forming a collinear central configuration. 
Unlike the convex kites, it is possible for three bodies of a concave kite central configuration to be nearly collinear.

Using MATLAB, we computed the index for each type of concave kite by evaluating $F(\wx, \wy)$ from the proof of Lemma~\ref{lemma:index} on a $300 \times 300$ grid in both
$\Cv_1$ and $\Cv_2$.  The results are shown in Figure~\ref{Fig:ConcaveDomain} and summarized in Table~\ref{table:Concave}.
Interestingly, there is a curve of degenerate central configurations (red curve) containing examples of each type.  These are configurations where $F = 0$ and
the modified Hessian is singular.

\renewcommand{\arraystretch}{2.2}
\begin{table}[h!]
\begin{center}
\arraycolsep=5pt
\begin{tabular}{|c|c|p{1.3in}|c|c|}
   \hline
   &  {\bf Mutual Distances}  & \hspace{0.4in} {\bf Index}    &  {\bf Degeneracies}   \\[0.05in]
   \hline  \hline
   {\bf Type 1}  &   $r_{12} > r_{24} > r_{23} > r_{13}$  &  \parbox{1.3in}{\hspace*{0.2in} $+1$ or $-1$ \\  \hspace*{0.1in} (majority $+1$)}   &  Yes   \\
   \hline
   {\bf Type 2}  &  $r_{24} > r_{12} > r_{13} > r_{23}$   &   \parbox{1.3in}{\hspace*{0.2in} $+1$ or $-1$ \\  \hspace*{0.1in} (majority $-1$)}   &  Yes   \\
   \hline
   {\bf $1+3$-gon}  & $r_{12} = r_{24}$ and $r_{23} = r_{13}$  & \parbox{1.2in}{$m > m^\ast \Longrightarrow +1$ \\  $m < m^\ast \Longrightarrow -1$} 
          &  $m^\ast = \frac{207 - 16\sqrt{3}}{338}$  \\[0.05in]
   \hline
\end{tabular}
\end{center}
\vspace{-0.15in}
\caption{Summary of information for the different types of concave kite central configurations.  We did not compute the index of the degenerate central configurations.}
\label{table:Concave}
\end{table}
\renewcommand{\arraystretch}{1.0}

Note that the degeneracy curve passes through the $1+3$-gon.  It is possible to compute the exact mass values for this particular degeneracy.  Plugging in
$\wx = \sqrt{3}$, $\wy = -1/\sqrt{3}$, $m_1 = m/2$, and $m_3 = 1 - 3m/2$ into inequality~(\ref{eq:condwMass}) yields a quadratic inequality in $m$
that is readily solved.  We find that the $1+3$-gon is degenerate for $m = m^\ast$, where
$$
m^\ast \; = \;  \frac{207 - 16\sqrt{3}}{338} \; \approx \;  0.530435 \, .
$$
From this we compute the ratio between the interior mass $m_3$ and one of the outer equal masses at the degeneracy to be $(81 + 64\sqrt{3})/249$, which matches
the value previously computed by Schmidt~\cite{schmidt} and Leandro~\cite{Leandro}. The $1+3$-gon is a local min for $m > m^\ast$ and a saddle for $m < m^\ast$.

\subsection{The mass map}

Next we consider the mass map $\wM$ from the configuration space $\Cv_1 \cup \, \Cv_2$  to the mass triangle $m_1 + m_3 + m = 1$, $m, m_1, m_3 \geq 0$.
Two coordinates of this map are given by
$$
m_1(\wx, \wy) \; = \;  \frac{\wy(s_{23} - s_{24}) (s_{23} - s_{13})}{S[ \wx(s_{13} - s_{12}) + \wy(s_{13} - s_{23}) ]}  \quad \mbox{and} \quad
m_3(\wx, \wy) \; = \;  \frac{\wx(s_{12} - s_{24}) (s_{12} - s_{13})}{S[ \wx(s_{13} - s_{12}) + \wy(s_{13} - s_{23}) ]} ,
$$
where $s_{ij} = r_{ij}^{-3}$ and $S = s_{13} + s_{24} - s_{12} - s_{23}$.  While each of these expressions are defined and positive on $\Cv_1 \cup \, \Cv_2$, they 
are undefined, taking the form $0/0$, at the $1+3$-gon $(x^\ast, y^\ast) = (\sqrt{3}, -1/\sqrt{3})$ because of its free mass parameter.
Nevertheless, limiting values of the masses exist at the $1+3$-gon.
Consider a line segment in $\Cv_1 \cup \, \Cv_2$ of slope $k$, with $k < - 2/3$, passing through $(x^\ast, y^\ast)$. Using a first-order Taylor expansion, it
is straight-forward to compute that 
$$
\lim_{(x,y) \rightarrow (x^\ast, y^\ast)} \; m_1(x,y) = \frac{6 + 9k}{18 - \sqrt{3} + 27k}  \quad \mbox{and} \quad
\lim_{(x,y) \rightarrow (x^\ast, y^\ast)} \; m_3(x,y) = \frac{-\sqrt{3}}{18 - \sqrt{3} + 27k} .
$$
The limiting values satisfy $m_3 = 1 - 3m_1$, which is equivalent to $m_1 = m/2$, as expected for the $1+3$-gon.

\begin{figure}[t]
\begin{center}
\includegraphics[height = 330bp]{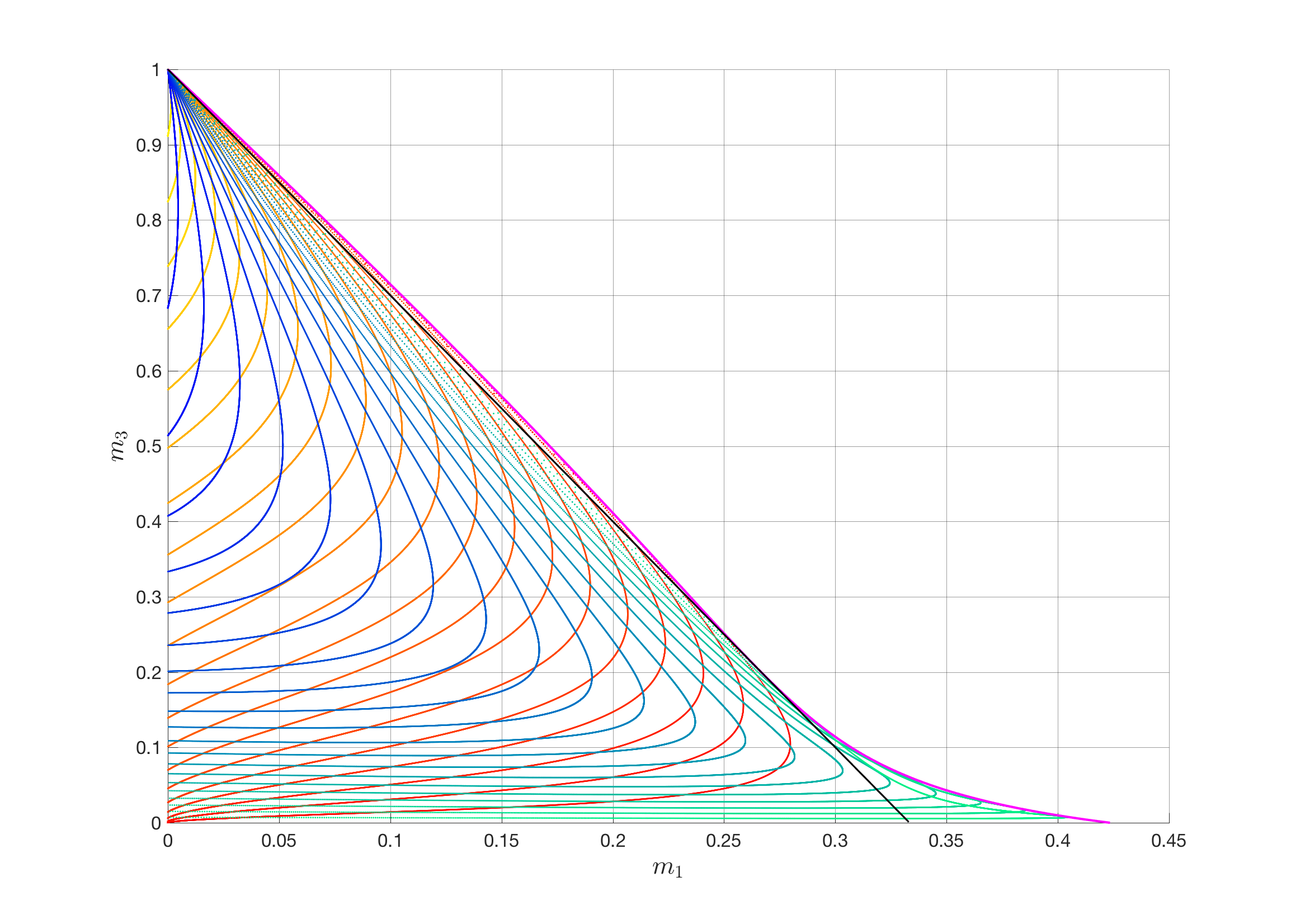}
\end{center}
\vspace{-0.45in}
\caption{The image of the different concave configurations under the mass map $\wM$ projected into the $m_1m_3$-plane.} 
\label{Fig:MassMap-m1m3}
\end{figure}

The image of the mass map $\wM$ projected into the $m_1m_3$-plane is shown in Figure~\ref{Fig:MassMap-m1m3}.  The figure was created by computing the image of 20
equally spaced vertical lines in both $\Cv_1$ (red and orange curves) and $\Cv_2$ (blue and green curves).  The image of the red degeneracy curve from Figure~\ref{Fig:ConcaveDomain}
is shown in magenta
and the masses for the $1+3$-gon are shown on the black line $m_3 = 1 - 3m_1$.  There are several interesting features revealed by this image.

\begin{figure}[t]
\begin{center}
\includegraphics[height = 300bp]{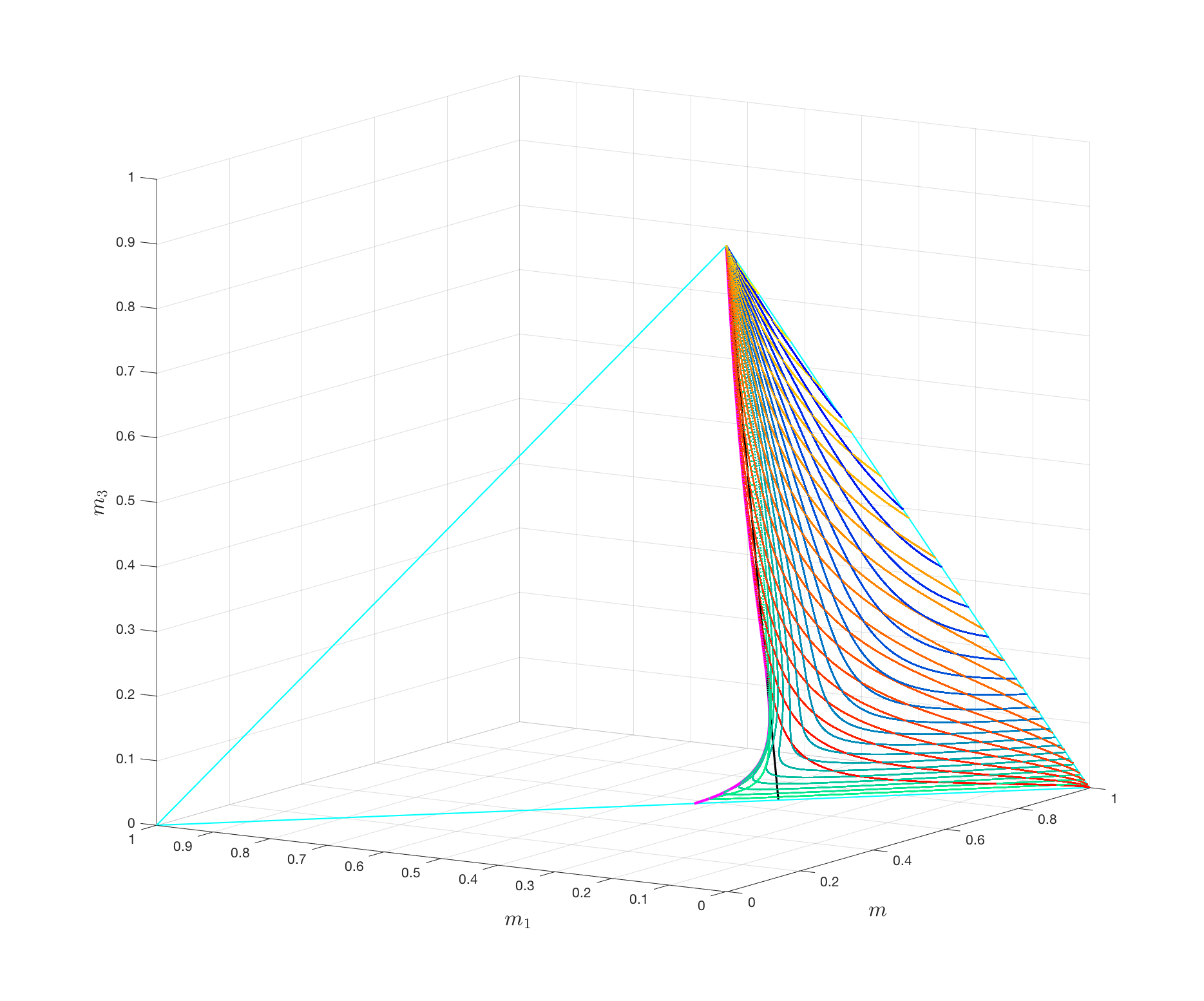}
\end{center}
\vspace{-0.45in}
\caption{The image of the mass map for the concave kites only covers a portion of the normalized mass triangle $m_1 + m_3 + m = 1$ (boundaries of this triangle shown in cyan).} 
\label{Fig:MassMap3D}
\end{figure}

\begin{itemize}
\item  $\wM$ does not cover the full mass space (see also Figure~\ref{Fig:MassMap3D}).  
There are mass values for which there does {\em not exist} a concave kite central configuration (e.g., $m_1 = 0.25, m_3 = 0.4$).  
This holds even if we interchange $m_1$ and $m_3$, and $\wx$ and $\wy$ (moving body 1 to the interior of the triangle).  
The value of $m_1$ is bounded above by $0.42345$ and for each $m_1$, there is a maximum allowable
value for the interior mass $m_3$ determined by the image of the degeneracy curve.

\item  As indicated by the intersections in the figure, $\wM$ is a two-to-one map.  For those mass values where a concave kite exists, 
most have one configuration of each type (e.g., a point of intersection between a blue curve and an orange curve). 
However, for mass values located above the black line and below the magenta curve, the two pre-images are of the {\em same} type.  This can be seen, for example, in the lower right
part of the figure, where different green curves intersect each other.  Moving across the black line in Figure~\ref{Fig:MassMap-m1m3}, 
one of the pre-images of $\wM$ will pass through the $1+3$-gon and then flip its type of configuration.

\item  Any curve that meets the magenta curve must be tangent to it.  This follows because the mass map is not invertible on the degeneracy curve.
The derivative of $\wM$ evaluated at any degenerate point is a rank one matrix. 
\end{itemize}

\section{Linear Stability}

In this section we numerically investigate the spectral and linear stability of the relative equilibria that arise 
from the kite central configurations.  We find that nearly all of the
kite relative equilibria are linearly unstable.  The only stable solutions occur for the convex case, very close to the lower boundary of $\C$.  Surprisingly, all of the concave solutions
are unstable, even those with a ``dominant'' mass (where one body is substantially larger than the others).  
In agreement with Moeckel's conjectures (see problems 15 and~16 in~\cite{albouy-pblms}), the linearly stable solutions found
correspond to nondegenerate local minima of $U|_{\N}$ and contain a dominant mass.  Specifically, we compute that the infimum of the mass ratio
$m_1/(m_2 + m_3 + m_4)$ over the set of linearly stable solutions is $(25 + 3\sqrt{69})/2 \approx 24.9599$.
Detailed plots showing the different eigenvalue structures for the two types of concave solutions are provided, revealing some interesting bifurcation curves.

\subsection{An invariant subspace for the nontrivial eigenvalues}

We first provide an overview and definition of linear stability of relative equilibria in the planar $n$-body problem (see~\cite{rick-lindom}).
Here we work in the full phase space, utilizing the Hamiltonian structure of the problem.
Given the correct initial momentum, a central configuration $q \in \mathbb{R}^{2n}$ will rigidly rotate about its center of mass with angular velocity $\omega = \sqrt{\lambda}$,
where $\lambda = U(q)/I(q)$ is the same value given in the defining c.c.~equation~(\ref{eq:cc-UI}).  After changing to rotating coordinates, 
the linear stability of this circular periodic solution is determined by regarding the relative equilibrium as a fixed point in the rotating frame.
As with our initial setup, we will assume that the center of mass of $q$
is at the origin and that $I(q) = 1$.  Linearizing about the relative equilibrium gives the $4n \times 4n$ {\em stability matrix}
$$
\Lambda \; = \;  
\begin{bmatrix}
\omega K   &  M^{-1} \\[0.04in]
D^2U(q)     &  \omega K 
\end{bmatrix}
,
$$
where $M = \mbox{diag }\{m_1, m_1, m_2, m_2, \ldots, m_n, m_n\}$ and $K$ is the $2n \times 2n$ block-diagonal matrix with
$J = 
\begin{bmatrix}
0 & 1\\
-1 & 0 
\end{bmatrix}
$
on the diagonal.  \\

The stability matrix has two ``trivial'' eigenspaces that arise from the symmetries of the problem.  The linear subspace
$$
W_1 \; = \;  \mbox{span } \{ (q,0), \, (Kq,0), \, (0, Mq), \, (0, KMq) \}
$$
is invariant under $\Lambda$ with eigenvalues $0, 0, \pm \omega i$.  These eigenvalues are a result of the rotational symmetry ($q$ is not isolated) and the homogeneity of~$U$.
The linear subspace
$$
W_2 \; = \;  \mbox{span } \{ (\xi_n,0), \, (K\xi_n,0), \, (0, M\xi_n), \, (0, KM\xi_n) \} ,
$$
where $\xi_n = (1, 0, 1, 0, \ldots, 1, 0)$, is also invariant under $\Lambda$ with eigenvalues $\pm \omega i, \pm \omega i$.  
These eigenvalues come from the conservation of the center of mass and linear momentum.

Thus, any relative equilibrium will contain the eight trivial eigenvalues $0,0, \pm \omega i, \pm \omega i, \pm \omega i$.  It is standard to call $q$ a {\em nondegenerate} 
relative equilibrium if the remaining $4n - 8$ eigenvalues are nonzero.  It is {\em spectrally stable} if the nontrivial eigenvalues lie on the imaginary axis.
Linear stability is defined using the skew inner product $\Omega$.  For vectors $v,w \in \mathbb{C}^{4n}$, let
$$
\Omega(v,w) \; = \;  v^{\footnotesize{T}} \mathbb{J} w, \; \; \mbox{where }
\mathbb{J} = 
\begin{bmatrix}
0 & I_{2n}\\
-I_{2n} & 0 
\end{bmatrix}
\mbox{ and }
I_{2n} \mbox{ is the } 2n \times 2n  \mbox{ identity matrix.}
$$
Let $W$ denote the skew-orthogonal complement of $W_1 \cup W_2$, that is,
$$
W \; = \;  \{v \in \mathbb{C}^{4n}: \Omega(v,w) = 0 \; \forall w \in W_1 \cup W_2 \} \, .
$$
Using the fact that $\Omega(v, \Lambda w) = -\Omega(\Lambda v, w)$, it is easy to see that $W$ is also an invariant subspace under~$\Lambda$.
The nontrivial eigenvectors all lie in $W$. 
Therefore, we call a relative equilibrium {\em linearly stable} if it is spectrally stable and if
$\Lambda$ is diagonalizable when restricted to~$W$.

We note that there exists a vector $\chi \in \mathbb{R}^{2n}$ in the null space of the 
modified Hessian $M^{-1} D^2U(q) + \lambda I_{2n}$ if and only if the vector $(\chi, -\omega KM \chi)$ is in the kernel of $\Lambda$.  
Hence $q$ is nondegenerate when regarded as a central configuration (in the same sense as Definition~\ref{def:CC})
if and only if it is a nondegenerate relative equilibrium.

Next we focus on the linear stability of the kite relative equilibria.  
In order to determine stability, we need to find the eight nontrivial eigenvalues of $\Lambda$.  Suppose that
$q = (a, 0, -c, d, -b, 0, -c, -d)$ is a kite central configuration with center of mass at the origin and moment of inertia $I(q) = 1$.  We will use the mass inner
product 
$$
<\!v, w\!> \;  =  \;  v^{\footnotesize{T}} M w, \quad \mbox{for } v,w \in \mathbb{R}^8 ,
$$
where $M = \mbox{diag }\{m_1, m_1, m_2, m_2, m_3, m_3, m_4, m_4\}$.  Recall that $m_2 = m_4 = m/2$ and that the masses have been normalized by
$m_1 + m_3 + m = 1$.  The vectors $q, Kq, \xi_4,$ and $K\xi_4$ are mutually orthogonal and of unit length with respect to this inner product.

More accurate numerical computations can be achieved by restricting $\Lambda$ to $W$, the skew orthogonal complement of $W_1 \cup W_2$, which
is an eight-dimensional vector space.  
In order to find a basis for $W$,
we first extend $q, Kq, \xi_4,$ and $K\xi_4$ to an $M$-orthonormal basis for $\mathbb{R}^8$.  
We take a geometric approach and use the oriented areas $\Delta_i$ of the triangles
obtained by deleting the $i$th body~\cite{rick-book}.  For our setup, we have
$$
\Delta_1 = -d(b-c), \; \Delta_3 = -d(a+c), \quad \mbox{and } \quad \Delta_2 = \Delta_4 = \casefrac{1}{2} d(a+b) .
$$
Define the vector $u_1 = (\Delta_1/m_1, 0, \Delta_2/m_2, 0, \Delta_3/m_3, 0, \Delta_4/m_4, 0)$. 
It is straight-forward to check that the vectors $u_1$ and $-K u_1$ are each $M$-orthogonal to $q, Kq, \xi_4,$ and $K\xi_4$.

Next we use Gram-Schmidt orthogonalization to obtain an orthogonal vector $u_2$.  Recall that 
$$
 \mbox{proj}_{w}(v) \; = \;  \frac{  <\!v, w\!> }{ <\!w, w\!> } w
 $$ 
is the projection of the vector $v$ onto $w$.  Suppose that $v \in \mathbb{R}^8$ is an arbitrary vector.  Define
\begin{equation}
u_2 \; = \;  v -  \left( \mbox{proj}_{q}(v)  + \mbox{proj}_{K q}(v)  + \mbox{proj}_{\xi_4}(v) + \mbox{proj}_{K \xi_4}(v) +  \mbox{proj}_{u_1}(v) + \mbox{proj}_{-K u_1}(v) \right)  .
\label{orthVecU2}
\end{equation}
The vector $u_2$ is $M$-orthogonal to the six vectors $q, Kq, \xi_4, K\xi_4, u_1,$ and $-K u_1$.
It is nonzero if $v$ does not lie in the span of these six vectors.

Taking $v = (1, 0, 0, 0, 0, 0, 0, 0)$ in~(\ref{orthVecU2}) yields the vector
\begin{eqnarray*}
u_2  & = &  (1 - m_1a^2 - m_1 - \frac{\Delta_1^2}{m_1  <\!u_1, u_1\!>}, \, 0,\, m_1ac - m_1 - \frac{\Delta_1 \Delta_2}{m_2  <\!u_1, u_1\!>}, \, -m_1 a d, \\[0.07in]
&  &  m_1ab - m_1 - \frac{\Delta_1 \Delta_3}{m_3  <\!u_1, u_1\!>}, \, 0, \, m_1ac - m_1 - \frac{\Delta_1 \Delta_4}{m_4  <\!u_1, u_1\!>}, \, m_1 a d) .
\end{eqnarray*}
Since $ad > 0$, we know that $u_2$ is nonzero.  
Finally, let $\widehat{u_1} = u_1/|| u_1 ||$ and $\widehat{u_2} = u_2/|| u_2 ||$ (norm taken with respect to our inner product). 
The vectors $\widehat{u_1}, -K\widehat{u_1}, \widehat{u_2},$ and $-K\widehat{u_2}$ form an $M$-orthonormal
basis for the orthogonal complement of span$\{ q, Kq, \xi_4, K\xi_4\}$ in $\mathbb{R}^8$.  
This in turn implies that the vectors
$$
(\widehat{u_1},0), \, (-K\widehat{u_1},0), \, (\widehat{u_2},0), \, (-K\widehat{u_2},0), \, (0, M \widehat{u_1}), \, (0, -MK \widehat{u_1}), \, (0, M \widehat{u_2}), \, (0, -MK \widehat{u_2}) 
$$
form a basis for $W$.
Using this basis, we see that the nontrivial eigenvalues are equal to the eigenvalues of the $8 \times 8$ matrix
$$
\Lambda_W \; = \;  
\begin{bmatrix}
\omega K & I_{4}\\[0.04in]
\widehat{A} & \omega K 
\end{bmatrix}
,
$$
where $\widehat{A}$ is the restriction of $M^{-1} D^2U(q)$ to span$\{\widehat{u_1}, -K\widehat{u_1}, \widehat{u_2}, -K\widehat{u_2}\}$.
The entries of $\widehat{A}$ can be computed numerically using our inner product.  For example, the first entry on the main diagonal
is simply $<\! \widehat{u_1}, M^{-1}D^2U(q) \widehat{u_1}\!>$.  Due to the symmetry of the kites and our choice of vectors, 
half the entries in $\widehat{A}$ vanish (specifically, $\widehat{a}_{ij} = 0$ if $i+j \equiv 1 \mod{2}$).

The matrix $\Lambda_W$ is Hamiltonian and therefore has an even characteristic polynomial.  The eigenvalues consist of
real pairs $\pm a_j$, pure imaginary pairs $\pm b_j \, i$, or complex quadruplets $\pm a_j \pm b_j \, i$.  
The relative equilibrium is linearly stable if it has four distinct imaginary pairs of eigenvalues.

\subsection{Convex case}

For the convex kites, most of the relative equilibria are unstable, except for a small neighborhood about the lower boundary in~$\C$.

\begin{prop}
There is an open neighborhood in ${\cal C}$ about the lower boundary curve $r_{12} = r_{13}$ for which the corresponding kite
relative equilibria are linearly stable. 
\end{prop}

\pf
This is a direct application of the work of Corbera et al.~(see Corollary~1 in~\cite{CCLM}).  The authors show that the convex critical points
of the $1+3$-body problem (which includes the lower boundary of $\C$) are nondegenerate minima of a special potential function $V(\theta)$.
By the main theorem in~\cite{rick-lindom}, this implies that for each point $(\wx, \wy = -\wx + \sqrt{\wx^2 + 1})$ on the lower boundary of~$\C$,
there is a smooth family of linearly stable relative equilibria converging to $(\wx, \wy)$.  These relative equilibria must be kites
because $m_2 = m_4$~\cite{albouy}.

Starting on the lower boundary of $\C$ and increasing $\wy$, stability can be lost in one of two ways.
Either a pair of pure imaginary eigenvalues collides at the origin and then splits into a real pair, 
or two pairs of pure imaginary eigenvalues merge on the imaginary axis and then break off into a complex quadruplet 
(a Krein collision).  In either case, we have repeated roots at the bifurcation and the discriminant of the characteristic polynomial of $\Lambda_W$
must vanish.  Since the entries of $\Lambda_W$ vary continuously in $(\wx, \wy)$, so to does this discriminant.  Thus there is an
open inequality which guarantees a linear stable family persists locally.
\enpf

\begin{figure}[h!]

\begin{center}
\includegraphics[width = 270bp]{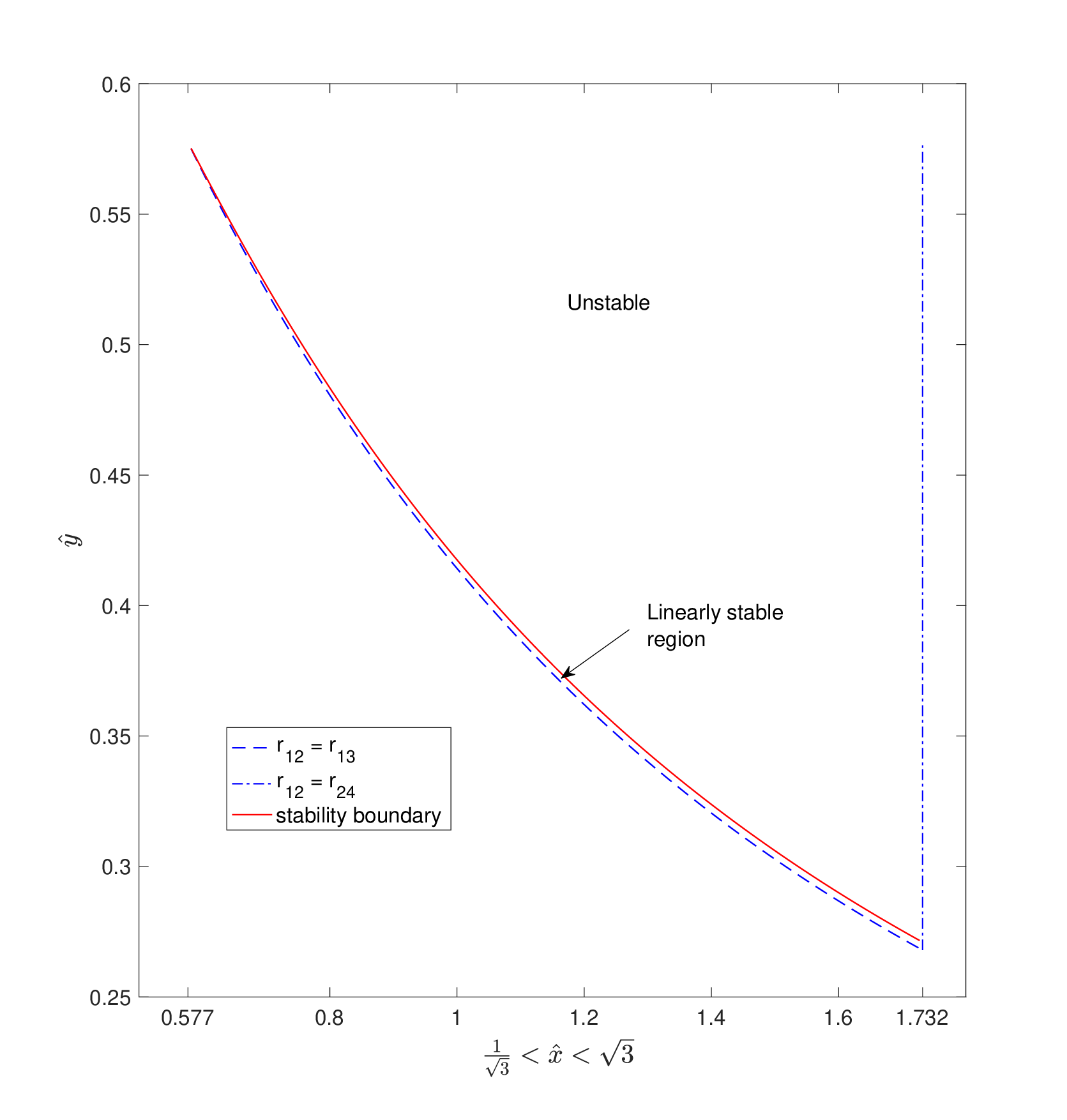}
\hspace{-0.25in}
\includegraphics[height = 200bp]{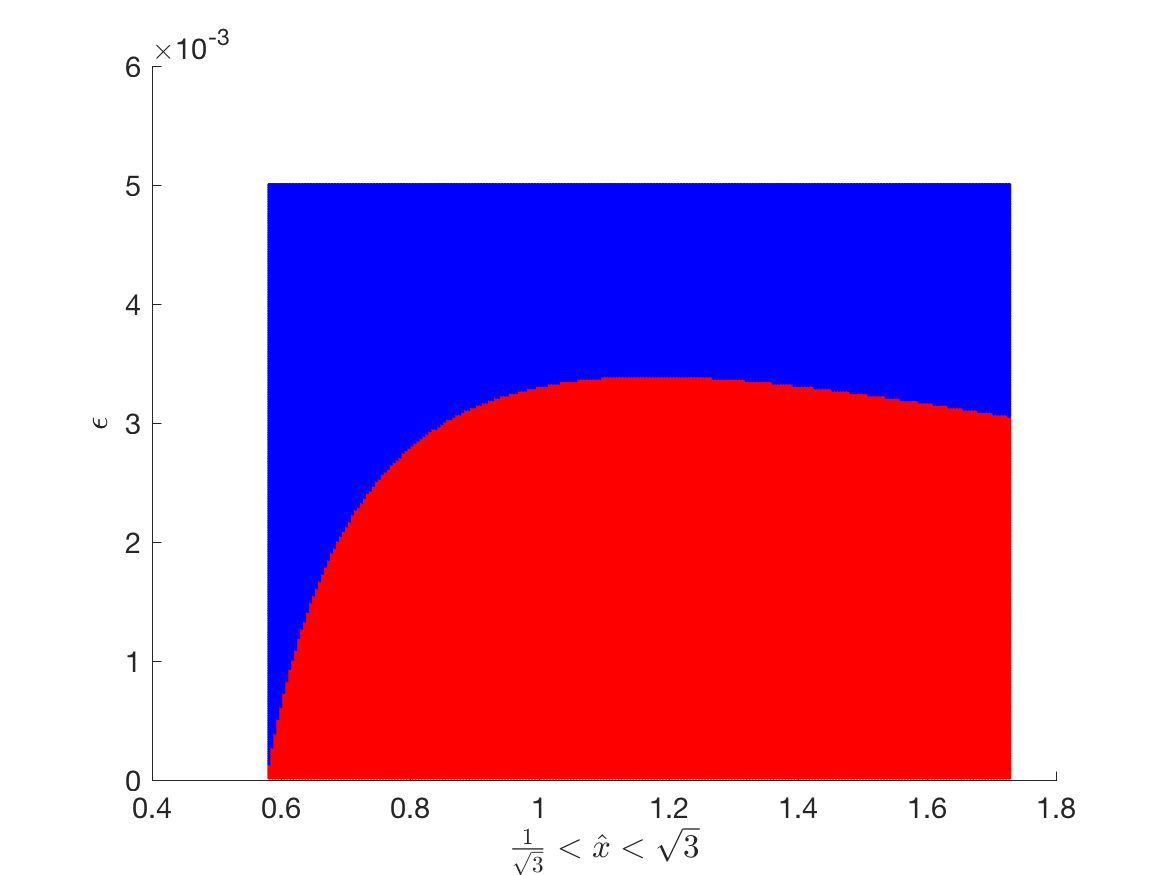}
\end{center}
\vspace{-0.2in}
\caption{Stability diagrams for the convex case, showing a small neighborhood of linearly stable relative equilibria near the lower boundary of $\C$. 
The right plot is a $250 \times 250$ grid of sample points within $0.005$ (vertically) of the lower boundary.  The vertical axis ($\epsilon$) indicates the
distance between $\hat{y}$ and the lower boundary.  Red points correspond to linearly stable solutions, while blue points are unstable.}
\label{Fig:StabPlotCvex}
\end{figure}

The stability diagram for the convex case is plotted in the reduced configuration space $\C$ in Figure~\ref{Fig:StabPlotCvex}.  The plots were created in MATLAB
using the \verb+eig+ command to compute the eigenvalues of $\Lambda_W$. Given a point $(\wx, \wy) \in \C$, we scale and shift the configuration
using equations (\ref{eq:a})--(\ref{eq:d}) to obtain the vector $q \in \mathbb{R}^8$ with center of mass at the origin and $I(q) = 1$.  The relative 
equilibrium is considered linearly stable if the real parts of the nontrivial eigenvalues are less than $1 \times 10^{-10}$ and the eigenvalues are distinct.
Eigenvalues were normalized by dividing by $\omega$.  Stability is lost through a Krein collision, which occurs on the red boundary curve shown in the left plot.  
As shown in the right plot, the size of the stable region is quite small, less than $0.004$ in the vertical direction.

\subsection{Concave case}

None of the concave kite relative equilibria are spectrally stable.  The configurations with a dominant mass are located close to the boundary
$r_{23} = r_{13}$, which contains critical points of the $1+3$-body problem (blue dashed curve in Figure~\ref{Fig:ConcaveDomain}).
The third body (largest mass) in the interior of the triangle is roughly equidistant to the other three smaller bodies.  
We checked numerically that the configurations on the boundary are either local maxima or saddles of the special potential function $V(\theta)$
defined in~\cite{rick-lindom}.  It follows by the main theorem in~\cite{rick-lindom} that the nearby kite relative equilibria (both types)
are linearly unstable.  Between the two types is the $1+3$-gon relative equilibrium 
(equilateral triangle with an arbitrary mass at the center), but this is known to be unstable for any value of the central mass~\cite{gr-ngon}.

\begin{figure}[h!]
\begin{center}
\includegraphics[width = 370bp]{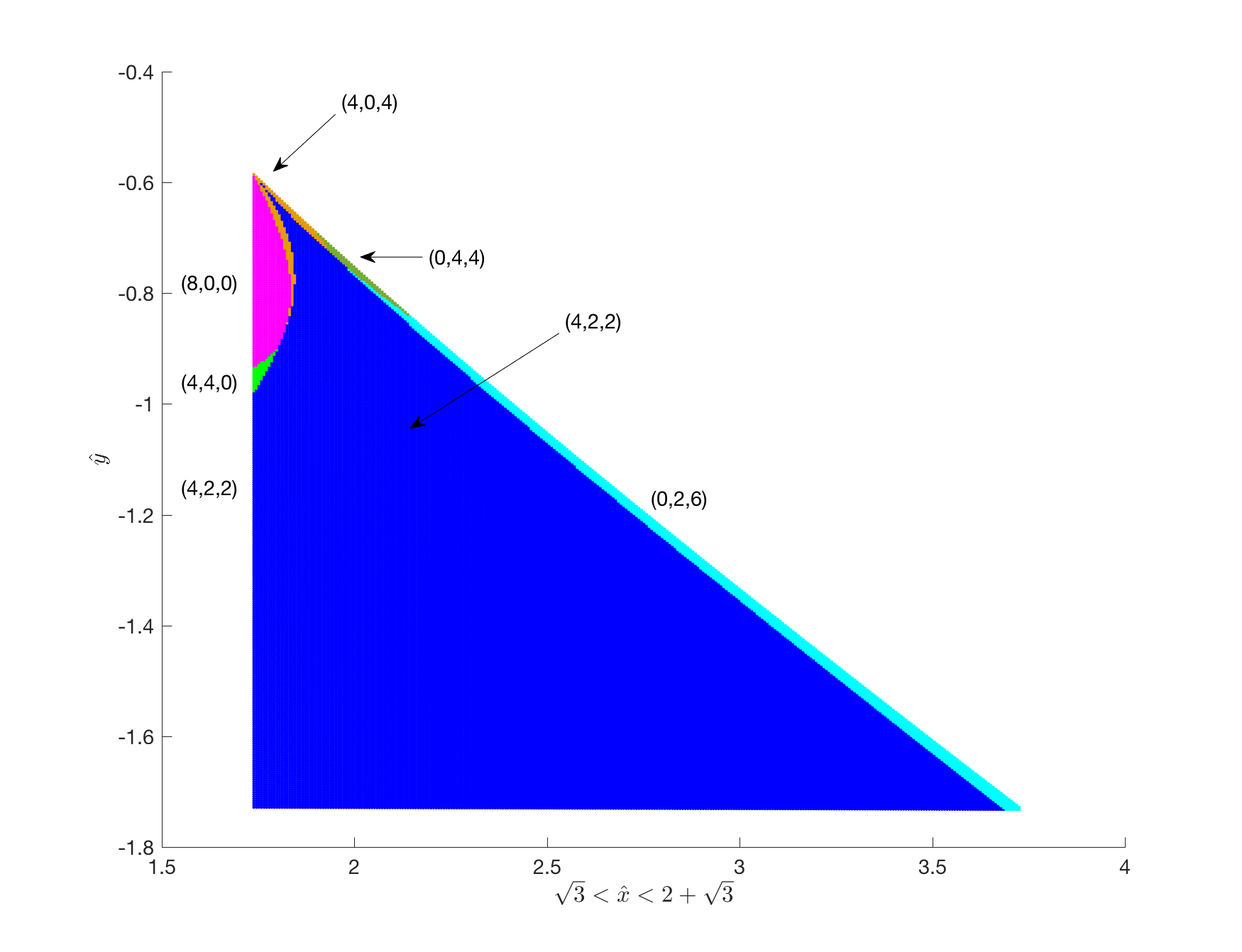}
\end{center}
\vspace{-0.4in}
\caption{Stability diagram for the Type 1 concave kite relative equilibria.  
Each colored region represents a different eigenvalue structure, as indicated by the triple 
(fully complex, real, pure imaginary).} 
\label{Fig:StabPlotCave1}
\end{figure}

\begin{figure}[t]
\begin{center}
\includegraphics[width = 390bp]{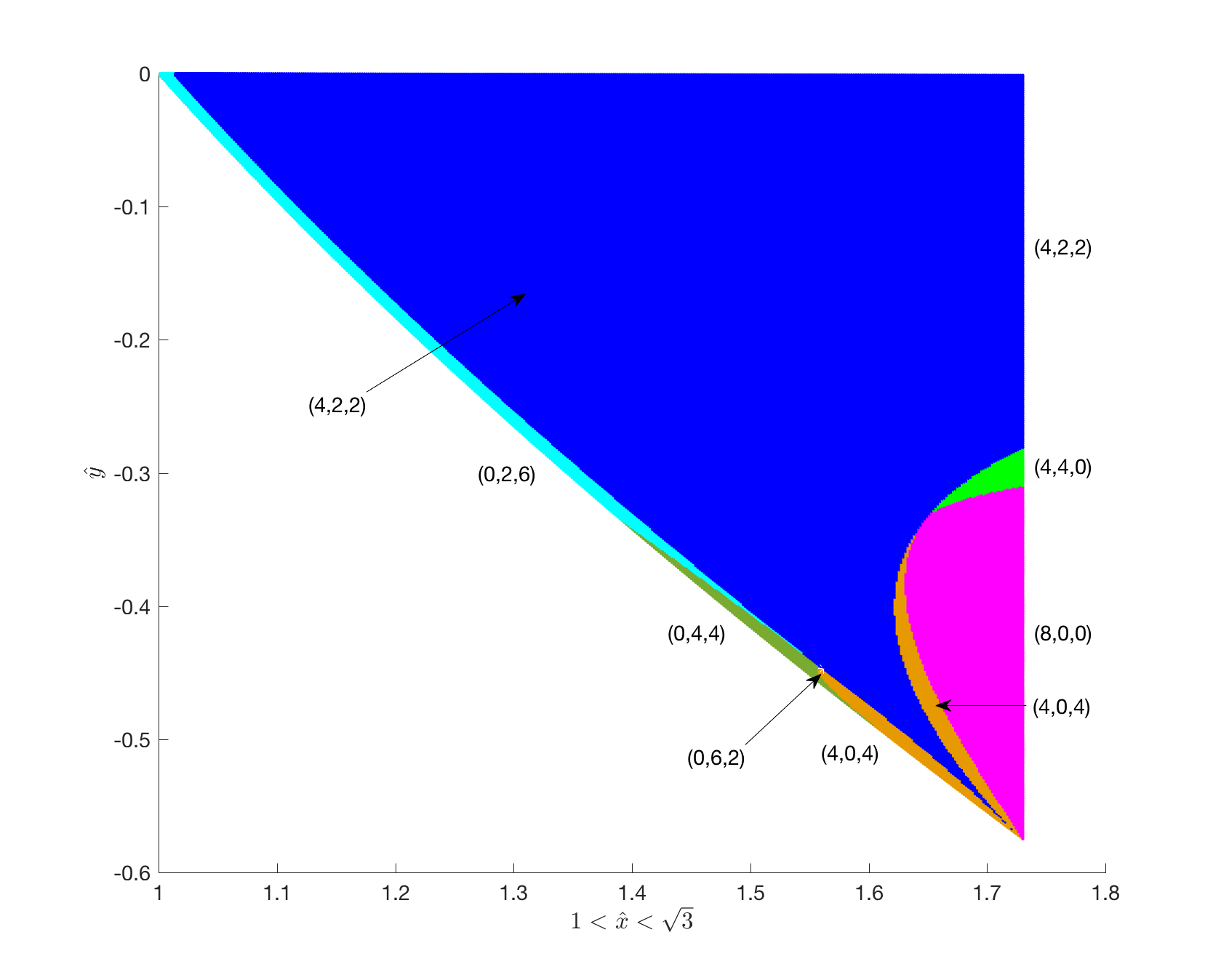}
\end{center}
\vspace{-0.4in}
\caption{Stability diagram for the Type 2 concave kite relative equilibria.} 
\label{Fig:StabPlotCave2}
\end{figure}

Detailed stability diagrams for both the Type 1 and 2 concave configurations were created using MATLAB
(see Figures \ref{Fig:StabPlotCave1} and~\ref{Fig:StabPlotCave2}).  
Each color represents a different structure for the eight nontrivial eigenvalues, as indicated by the triple 
(fully complex, real, pure imaginary).  For example, eigenvalues for points in the magenta regions $(8,0,0)$ consist of two complex quadruplets, 
while the eigenvalues in the large blue regions $(4,2,2)$ consist of a complex quadruplet, a real pair, and a pure imaginary pair.
Two of the nine possible combinations of eigenvalues do not exist for either type: four real pairs $(0,8,0)$ and the stable
case of four pure imaginary pairs $(0,0,8)$.

There are several interesting bifurcation curves apparent in Figures \ref{Fig:StabPlotCave1} and~\ref{Fig:StabPlotCave2}, 
where the eigenvalue structure changes.  Some of these 
curves match the degenerate curves shown in Figure~\ref{Fig:ConcaveDomain}. For example, moving from an orange region $(4,0,4)$ into a blue region $(4,2,2)$, 
one pair of pure imaginary eigenvalues merge at the origin (a degeneracy) and then switch over to a pair of real eigenvalues.  
Krein collisions occur when passing from a cyan region $(0,2,6)$ into a blue region $(4,2,2)$, or traveling from an orange region $(4,0,4)$ 
into a magenta region $(8,0,0)$.
A third example occurs when two pairs of real eigenvalues merge on the real axis and then break off into a complex quadruplet,
a scenario that transpires when moving from a neon green region $(4,4,0)$ into a magenta region $(8,0,0)$.  For those bifurcations 
where the eigenvalues pass through the origin (a degeneracy), but which are not on the curve shown in Figure~\ref{Fig:ConcaveDomain}, 
these must feature eigenvectors that break the symmetry of the kite (perturbations that exit $\Cv_1$ or $\Cv_2$).
The corresponding kite central configurations are degenerate in the full four-body problem.  

\subsection{Relative size of the dominant mass}

Our final computation explores the masses of the linearly stable solutions, focusing on the size of the largest mass with respect to the sum of the others.
For the stable convex kites, the largest mass is $m_1$.  Consider the mass ratio 
$$
\psi \; = \; \frac{m_1}{m_2 + m_3 + m_4} \; = \;  \frac{m_1}{1 - m_1} \, 
$$
evaluated on the set of all linearly stable solutions.
This ratio is infinite on the lower boundary of $\C$ since $m = m_3 = 0$ there, and numerical experiments show that it is decreasing as a function of $\wy$.  
It follows that the minimum value of $\psi$ will occur on the red boundary curve from the left plot in Figure~\ref{Fig:StabPlotCvex}.

\begin{figure}[t]

\begin{center}
\includegraphics[height = 280bp]{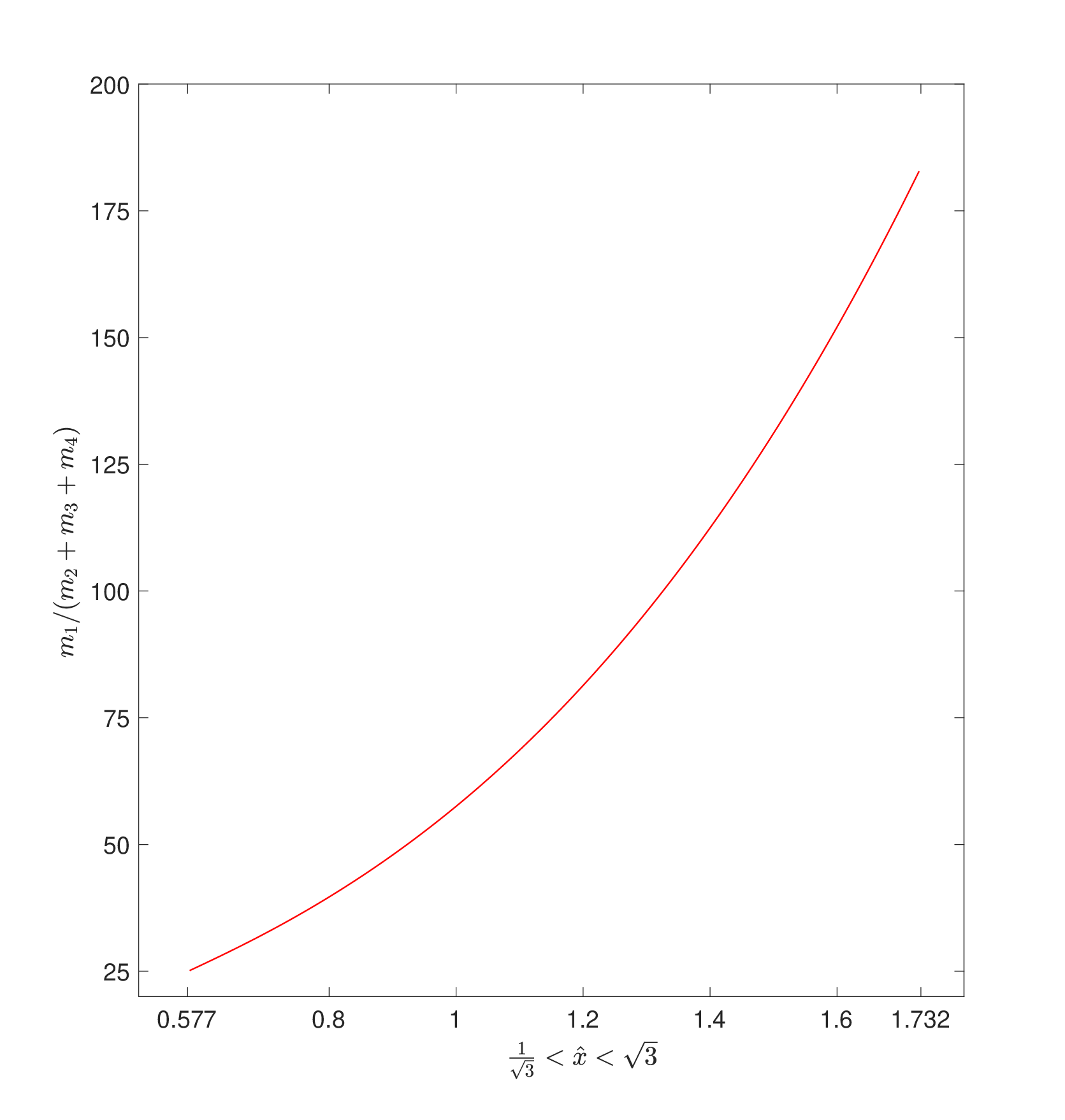}
\end{center}
\vspace{-0.3in}
\caption{The mass ratio $\psi = m_1/(m_2 + m_3 + m_4)$ evaluated along the stability boundary curve from the left plot in Figure~\ref{Fig:StabPlotCvex}.
For each $\wx$, in order for the convex kite relative equilibrium to be linearly stable, the mass ratio must be greater than the value given by this curve.}
\label{Fig:MassRatCvex}
\end{figure}

A graph of $\psi$ evaluated along the stability boundary is shown in Figure~\ref{Fig:MassRatCvex}.  It is an increasing function of $\wx$,
with a minimum value around 25.  Careful numerical calculations near the boundary point $(\wx, \wy) = (1/\sqrt{3}, 1/\sqrt{3} \,)$ indicate that the
infimum of $\psi$ is 
\begin{equation}
\frac{25 + 3 \sqrt{69}}{2} \; \approx \; 24.959936 \, .
\label{inf-massrat}
\end{equation}
The exact value is obtained from Routh's critical mass ratio $\rho_r = (1 - \sqrt{69}/9)/2$ for the restricted three-body problem,
where the two primaries have masses $\rho$ and $1 - \rho$. Placing an infinitesimal mass at each of the Lagrange triangle points $L_4$ and $L_5$ 
yields a four-body central configuration that is equivalent to the kite c.c.~for $(\wx, \wy) =  (1/\sqrt{3}, 1/\sqrt{3} \,)$. 
This is a limiting configuration for a linearly stable family of relative equilibria in the four-body problem, provided $\rho < \rho_r$~\cite{rick-clusters, gr-thesis}.
Substituting $m_1 = 1 - \rho_r$ into the mass ratio $\psi$ gives the value shown in~(\ref{inf-massrat}).
These numerical results provide further evidence for Moeckel's dominant mass conjecture.

\vs

\nin {\bf Acknowledgements:}  
The author would like to thank Dave Damiano, Victor Donnay, John Little, and Manuele Santoprete for their support and suggestions regarding this work. 
The value of the mass ratio~$\psi$ discussed in the last section was first investigated by undergraduate researchers 
Margaret Hauser and Gopal Yalla during a summer research project supervised by the author, with financial support form NSF grant DMS-1211675.

\bibliographystyle{amsplain}

\end{document}